\documentclass[12pt,a4paper,twoside]{article}

\usepackage[intlimits]{amsmath}
\usepackage{amssymb}
\usepackage{fancyhdr}

\usepackage{a4wide}

\newcommand{\maple}{Maple }
\newcommand{\mapleSE}{Maple}

\newtheorem{theorem}{Theorem}
\newtheorem{definition}[theorem]{Definition}
\newtheorem{example}[theorem]{Example}
\newtheorem{remark}[theorem]{Remark}

\newenvironment{keywords}{\begin{center}
\begin{minipage}[c]{130mm} {\bf Keywords:}} {\end{minipage}
\end{center}}
\newenvironment{msc}{\begin{center}
\begin{minipage}[c]{130mm} {\bf 2000 Mathematics Subject Classification:}} {\end{minipage}
\end{center} \bigskip}

\begin{document}

\title{Automatic Computation of Conservation Laws\\
in the Calculus of Variations and Optimal Control\thanks{Research report CM05/I-24,
Dept. Mathematics, Univ. Aveiro, June 2005. Partially
presented at the 10th International Conference \emph{Mathematical
Modelling and Analysis}, and 2nd International Conference
\emph{Computational Methods in Applied Mathematics}, June 1 - 5,
2005, Trakai, Lithuania. Submitted to Journal
``Computational Methods in Applied Mathematics''. A Maple worksheet
is available at \texttt{http://www.mat.ua.pt/delfim/maple.htm}}}

\author{Paulo~D.~F.~Gouveia\\ \texttt{pgouveia@ipb.pt} \and
Delfim~F.~M.~Torres\\ \texttt{delfim@mat.ua.pt}}

\date{Control Theory Group (\textsf{cotg})\\
Department of Mathematics\\
University of Aveiro\\
3810-193 Aveiro, Portugal}

\maketitle


\begin{abstract}
We present analytic computational tools that permit us to
identify, in an automatic way, conservation laws in optimal
control. The central result we use is the famous Noether's
theorem, a classical theory developed by Emmy Noether in 1918, in
the context of the calculus of variations and mathematical
physics, and which was extended recently to the more general
context of optimal control. We show how a Computer Algebra System
can be very helpful in finding the symmetries and corresponding
conservation laws in optimal control theory, thus making useful in
practice the theoretical results recently obtained in the
literature. A Maple implementation is provided and several
illustrative examples given.
\end{abstract}



\begin{msc}
49K15; 49-04; 49S05.
\end{msc}


\section{Introduction}

Optimal control problems are usually solved with the help of the
famous Pontryagin maximum principle \cite{Pontryagin62}, which is
a generalization of the classic Euler-Lagrange and Weierstrass
necessary optimality conditions of the calculus of variations. The
method of finding optimal solutions via Pontryagin's maximum
principle proceeds through the following main three steps: (i) one
defines the Hamiltonian of the problem; (ii) with the help of the
maximality condition one tries to express the control variables
with respect to state and adjoint variables; (iii) the Hamiltonian
system is written in terms of state and adjoint variables only,
and the solutions of this system of ordinary differential
equations are sought. Steps (ii) and (iii) are, generally
speaking, nontrivial, and very difficult (or even impossible) to
implement in practice \cite{Serovaiskii}. One way to address the
problem is to find conservation laws, \textrm{i.e.}, quantities
which are preserved along the extremals of the problem. Such
conservation laws can be used to simplify the problem
\cite{GrizzleMarcus,alik}. The question is then the following: how
to determine these conservation laws? It turns out that the
classic results of Emmy Noether \cite{Noether,OLV} of the calculus
of variations, relating the existence of conservation laws with
the existence of symmetries, can be generalized to the wider
context of optimal control \cite{Djukic73,Gogodze88,delfimEJC},
reducing the problem to the one of discovering the
invariance-symmetries. The difficulty resides precisely in the
determination of the variational symmetries. While in Physics and
Economics the question of existence of conservation laws is
treated in a rather natural way, because the application itself
suggest the symmetries (\textrm{e.g.}, conservation of energy,
conservation of momentum, income/health law, etc -- all of them
coming from very intuitive symmetries of the problem), from a
strictly mathematical point of view, given a problem of optimal
control, it is not obvious and not intuitive how one might derive
a conservation law. Therefore, it would be of great practical use
to have at our disposal computational means for the automatic
identification of the symmetries of the optimal control problems
\cite{GrizzleMarcus,alik}. This is the motivation of the present
work: to present a \maple package that can assist in this respect.
The results extend the previous investigations of the authors,
done in the classical context of the calculus of variations
\cite{gouv04}, to the more general and interesting setting of
optimal control \cite{Pontryagin62}, where the application of
symmetry and conservation laws is an area of current research
\cite{alik,3ncnw}.

The use of symbolic mathematical software is becoming, in recent
years, an effective tool in mathematics \cite{Richards}. Computer
algebra, also known as symbolic computation, is an
interdisciplinary area of mathematics and computer science.
Computer Algebra Systems, \maple as an example, facilitate the
interplay of conventional mathematics with computers. They are, in
some sense, changing the way we learn, teach, and do research in
mathematics \cite{BaileyBorwein}. They can perform a myriad of
symbolic mathematical operations, like analytic differentiation,
integration of algebraic formulae, factoring polynomials,
computing the complex roots of analytic functions, computing
Taylor series expansions of functions, finding analytic solutions
of ordinary or partial differential equations, etc. It is not a
surprise that they are becoming popular in control theory and
control engineering applications \cite{IEEE05}. Here we use the
\maple 9.5 system to find symmetries and corresponding
conservation laws in optimal control. The paper is organized as
follows. In \S\ref{sec:CL} the problem of optimal control is
introduced, the necessary definitions are given, and Noether's
theorem is introduced. The method of computing conservation laws
in optimal control is explained in \S\ref{sec:CCL}, and the \maple
package is then applied in \S\ref{sec:EI} for computing symmetries
and families of conservations laws to a diverse range of optimal
control problems. In \S\ref{sec:CV} we focus attention to the
symbolic computation of conservation laws in the calculus of
variations. The \maple procedures are given in
\S\ref{sec:ProcMaple}, and we end \S\ref{concl} with some comments
and directions of future work.


\section{Conservation Laws in Optimal Control}
\label{sec:CL}

The optimal control problem consists in the minimization of an
integral functional,
\begin{equation}
\label{eq:funcionalCO}
I[\mathbf{x}(\cdot),\mathbf{u}(\cdot)]
= \int_{a}^{b} L(t,\mathbf{x}(t),\mathbf{u}(t)) \,\textrm{d}t \, ,
\end{equation}
subject to a control system described by ordinary differential
equations,
\begin{equation}
\label{eq:sistCont}
\dot{\mathbf{x}}(t)= \boldsymbol{\varphi}(t,\mathbf{x}(t),\mathbf{u}(t)) \, ,
\end{equation}
together with appropriate boundary conditions. The Lagrangian
$L(\cdot,\cdot,\cdot)$ is a real function, assumed to be
continuously differentiable in $[a, b]
\times\mathbb{R}^n\times\mathbb{R}^m$; $t\in\mathbb{R}$ the
independent variable; $\mathbf{x}\!:[a, b]
\rightarrow\mathbb{R}^n$ the vector of state variables;
$\mathbf{u}\!:[a, b]\rightarrow\Omega\subseteq\mathbb{R}^m$,
$\Omega$ an open set, the vector of controls, assumed to be
piecewise continuous functions; and $\boldsymbol{\varphi}\!:[a, b]
\times\mathbb{R}^n\times\mathbb{R}^m \rightarrow\mathbb{R}^n$ the
velocity vector, assumed to be a continuously differentiable
vector function. We propose a computational method that permits to
obtain conservation laws for a given optimal control problem. Our
method is based in the version of Noether's theorem established in
\cite{Djukic73} (see also \cite{Torres05}). To describe a
systematic method to compute conservation laws, first we need to
recall the standard definitions of \emph{extremal}, and
\emph{conservation law}.

The central result in optimal control theory is the famous
Pontryagin maximum principle \cite{Pontryagin62}, which gives a
necessary optimality condition for the problems of optimal
control.
\begin{theorem}[Pontryagin maximum principle]
\label{thm:pMaxPont} If
$\left(\mathbf{x}(\cdot),\mathbf{u}(\cdot)\right)$ is a solution
of the optimal control problem
(\ref{eq:funcionalCO})-(\ref{eq:sistCont}), then there exists a
non-zero pair $\left(\psi_0,\boldsymbol{\psi}(\cdot)\right)$,
where $\psi_0 \leq 0$ is a constant and $\boldsymbol{\psi}(\cdot)$
an $n$-vectorial piecewise $C^1$-smooth function with domain $[a,
b]$, in such a way the quadruple
$\left(\mathbf{x}(\cdot),\mathbf{u}(\cdot),\psi_0,
\boldsymbol{\psi}(\cdot)\right)$ satisfies the following
conditions in almost all points $t$ in the interval $[a, b]$:
\renewcommand{\theenumi}{\roman{enumi}}
\renewcommand{\labelenumi}{\emph{(\theenumi)}}
\begin{enumerate}
\item the Hamiltonian system
\begin{eqnarray}
\label{eq:sistCont2}
\dot{\mathbf{x}}(t)\emph{\textsuperscript{T}}&=& \frac{\partial
H}{\partial
\boldsymbol{\psi}}(t,\mathbf{x}(t),\mathbf{u}(t),\psi_0,
\boldsymbol{\psi}(t))\, ,\\
\label{eq:sistAdj}
\dot{\boldsymbol{\psi}}(t)\emph{\textsuperscript{T}}&=&-
\frac{\partial H}{\partial
\mathbf{x}}(t,\mathbf{x}(t),\mathbf{u}(t),\psi_0,
\boldsymbol{\psi}(t)) \, ,
\end{eqnarray}
\item the maximality condition
\begin{equation}
\label{eq:condMax}
H(t,\mathbf{x}(t),\mathbf{u}(t),\psi_0,\boldsymbol{\psi}(t))=
{\mathop {\max}\limits_{\mathbf{v} \in \Omega}} \,
H(t,\mathbf{x}(t),\mathbf{v},\psi_0,\boldsymbol{\psi}(t)),
\end{equation}
\end{enumerate}
with Hamiltonian
\begin{equation}
\label{eq:hamilt}
H(t,\mathbf{x},\mathbf{u},\psi_0,\boldsymbol{\psi})= \psi_0
L(t,\mathbf{x},\mathbf{u})+
\boldsymbol{\psi}\emph{\textsuperscript{T}} \cdot
\boldsymbol{\varphi}(t,\mathbf{x},\mathbf{u}) \, .
\end{equation}
\end{theorem}

\begin{remark}
The right-hand side of equations (\ref{eq:sistCont2}) and
(\ref{eq:sistAdj}), of the Hamiltonian system, represent a
row-vector formed by the partial derivatives of the Hamiltonian
scalar function $H$ with respect to each one of the components of
the derivation variable. Equation (\ref{eq:sistCont2}) is nothing
more than the control system (\ref{eq:sistCont}); equation
(\ref{eq:sistAdj}) is known as the \emph{adjoint system}.
\end{remark}

\begin{definition}
A quadruple $\left(\mathbf{x}(\cdot),\mathbf{u}(\cdot),\psi_0,
\boldsymbol{\psi}(\cdot)\right)$ satisfying the Pontryagin maximum
principle is said to be a (Pontryagin) \emph{extremal}. An
extremal is \emph{normal} when $\psi_0\neq0$, \emph{abnormal} if
$\psi_0=0$.
\end{definition}

\begin{remark}
Since we assume $\Omega$ to be an open set, the maximality
condition \eqref{eq:condMax} implies the \emph{stationary
condition}
\begin{equation}
\label{eq:stat:cond}
\frac{\partial H}{\partial \mathbf{u}}(t,\mathbf{x}(t),\mathbf{u}(t),\psi_0,\boldsymbol{\psi}(t))
= \mathbf{0} \, , \quad t \in [a,b] \, .
\end{equation}
Using the Hamiltonian system \eqref{eq:sistCont2}-\eqref{eq:sistAdj},
and the stationary condition \eqref{eq:stat:cond},
it follows that, along the extremals, the total derivative
of the Hamiltonian with respect to $t$ equals its partial derivative,
$t \in [a,b]$:
\begin{equation}
\label{eq:prop:dHdt}
\frac{\mathrm{d}}{\mathrm{d}t} H(t,\mathbf{x}(t),\mathbf{u}(t),\psi_0,\boldsymbol{\psi}(t)) =
\frac{\partial H}{\partial t}(t,\mathbf{x}(t),\mathbf{u}(t),\psi_0,\boldsymbol{\psi}(t)) \, .
\end{equation}
\end{remark}

Let us now consider a one-parameter group of $\mathbb{C}^1$
transformations $\mathbf{h}^s:[a, b] \times \mathbb{R}^n \times
\mathbb{R}^m \times \mathbb{R} \times \mathbb{R}^n \rightarrow
\mathbb{R} \times \mathbb{R}^n \times \mathbb{R}^m \times
\mathbb{R}^n$ of the form
\begin{equation}
\label{eq:transf}
\mathbf{h}^s(t,\mathbf{x},\mathbf{u},\psi_0,\boldsymbol{\psi})=
(h_t^s(t,\mathbf{x},\psi_0,\boldsymbol{\psi}),
\mathbf{h}_\mathbf{x}^s(t,\mathbf{x},\psi_0,\boldsymbol{\psi}),
\mathbf{h}_\mathbf{u}^s(t,\mathbf{x},\mathbf{u},\psi_0,\boldsymbol{\psi}),
\mathbf{h}_{\boldsymbol{\psi}}^s(t,\mathbf{x},\mathbf{u},
\psi_0,\boldsymbol{\psi})) \, ,
\end{equation}
which reduces to the identity transformation
when the parameter $s$ vanishes:
\[
h_t^0(t,\mathbf{x},\psi_0,\boldsymbol{\psi})=t,\;
\mathbf{h}_\mathbf{x}^0(t,\mathbf{x},\psi_0,\boldsymbol{\psi})=
\mathbf{x},\;
\mathbf{h}_\mathbf{u}^0(t,\mathbf{x},\mathbf{u},\psi_0,\boldsymbol{\psi})=
\mathbf{u},\;
\mathbf{h}_{\boldsymbol{\psi}}^0(t,\mathbf{x},\mathbf{u},\psi_0,
\boldsymbol{\psi})=\boldsymbol{\psi} \, .
\]
Associated with a one-parameter group of transformations
(\ref{eq:transf}), we introduce the infinitesimal \emph{generators}
\begin{eqnarray}
\label{eq:transf2}
T(t,\mathbf{x},\psi_0,\boldsymbol{\psi})  =
\left. \frac{\partial} {\partial{s}} h^s_t(t,\mathbf{x},\psi_0\,
\boldsymbol{\psi})\right|_{s=0}\textrm{, }
\mathbf{X}(t,\mathbf{x},\psi_0,\boldsymbol{\psi})  = \left.
\frac{\partial}{\partial{s}} \mathbf{h}_\mathbf{x}^s(t,\mathbf{x},
\psi_0,\boldsymbol{\psi})\right|_{s=0}\textrm{,}\nonumber\\
\mathbf{U}(t,\mathbf{x},\mathbf{u},\psi_0,\boldsymbol{\psi})  =
\left. \frac{\partial} {\partial{s}} \mathbf{h}_\mathbf{u}^s
(t,\mathbf{x},\mathbf{u},\psi_0,\boldsymbol{\psi})\right|_{s=0}
\textrm{, }
\boldsymbol{\Psi}(t,\mathbf{x},\mathbf{u},\psi_0,\boldsymbol{\psi})
=  \left. \frac{\partial}{\partial{s}}
\mathbf{h}_{\boldsymbol{\psi}}^s(t,\mathbf{x},
\mathbf{u},\psi_0,\boldsymbol{\psi})\right|_{s=0}
\end{eqnarray}
Emmy Noether was the first who established a relation between the
existence of invariance transformations of the problem and the
existence of conservation laws \cite{Noether}. This relation
constitutes a universal principle that can be formulated, as a
theorem, in several different contexts, under several different
hypotheses (see \textrm{e.g.}
\cite{Djukic73,Gogodze88,alik,LOG77,OLV,TOR,Torres05}).
Contributions in the literature go, however, further than
extending Noether's theorem to different contexts, and weakening
its assumptions. Since the pioneering work by Noether
\cite{Noether}, several definitions of invariance have been
introduced for the problems of the calculus of variations (see
\textrm{e.g.} \cite{LOG77,LOG,MR83c:70020,TOR}); and for the
problems of optimal control (see \textrm{e.g.}
\cite{Djukic73,Gogodze88,delfimEJC,Torres04}). All these
definitions are given with respect to a one-parameter group of
transformations \eqref{eq:transf}. Although written in a different
way (some of these invariance/symmetry notions involve the
integral functional, others only the integrand; some of them
involve the original problem and the transformed one, others only
the rate of change with respect to the parameter; etc) it turns
out that, when written in terms of the generators
\eqref{eq:transf2}, one gets a necessary and sufficient condition
of invariance that, essentially, coincides with all those
definitions. For this reason, here we define invariance directly
in terms of the generators \eqref{eq:transf2}.
\begin{definition}[\cite{Djukic73,Torres05}]
We say that an optimal control problem
(\ref{eq:funcionalCO})-(\ref{eq:sistCont}) is \emph{invariant}
under \eqref{eq:transf2} or, equivalently, that \eqref{eq:transf2}
is a \emph{symmetry} of the problem, if, and only if,
\begin{equation}
\label{eq:detGerad} \frac{\partial H}{\partial t}T +\frac{\partial
H}{\partial \mathbf{x}}\cdot \mathbf{X} +\frac{\partial
H}{\partial \mathbf{u}}\cdot \mathbf{U} +\frac{\partial
H}{\partial \boldsymbol{\psi}}\cdot \boldsymbol{\Psi}
-\boldsymbol{\Psi}\textsuperscript{T} \cdot \dot{\mathbf{x}}
-\boldsymbol{\psi}\textsuperscript{T} \cdot \frac{\mathrm{d}
\mathbf{X}}{\mathrm{d}t} +H \frac{\mathrm{d}T}{\mathrm{d}t}=0\, ,
\end{equation}
with $H$ the Hamiltonian \eqref{eq:hamilt}.
\end{definition}
A symmetry is an intrinsic property of the optimal control problem
(\ref{eq:funcionalCO})-(\ref{eq:sistCont}) (an intrinsic property
of the corresponding Hamiltonian \eqref{eq:hamilt}), and does not
depend on the extremals. If one restricts attention to the
quadruples $(\mathbf{x}(\cdot),\mathbf{u}(\cdot),\psi_0,
\boldsymbol{\psi}(\cdot))$ that satisfy the Hamiltonian system and
the maximality condition, one arrives at E.~Noether's theorem:
along the extremals, equalities \eqref{eq:sistCont2},
\eqref{eq:sistAdj}, \eqref{eq:stat:cond}, and \eqref{eq:prop:dHdt}
permit to simplify \eqref{eq:detGerad} to the form
\begin{equation}
\label{eq:preLC} \frac{\mathrm{d}H}{\mathrm{d}t}T -
\dot{\boldsymbol{\psi}}\textsuperscript{T} \cdot \mathbf{X} -
\boldsymbol{\psi}\textsuperscript{T} \cdot \frac{\mathrm{d}
\mathbf{X}}{\mathrm{d}t} + H \frac{\mathrm{d}T}{\mathrm{d}t}=0
\Leftrightarrow \frac{\mathrm{d}}{\mathrm{d}t} \left( H T -
\boldsymbol{\psi}\emph{\textsuperscript{T}} \cdot \mathbf{X}
\right) = 0 \, .
\end{equation}
\begin{definition}
A function $C(t,\mathbf{x},\mathbf{u},\psi_0,\boldsymbol{\psi})$
which is preserved along all the extremals of the optimal control problem
and all $t \in [a,b]$,
\begin{equation}
\label{eq:def:CL}
C(t,\mathbf{x}(t),\mathbf{u}(t),\psi_0,\boldsymbol{\psi}(t)) =
\text{const}\, ,
\end{equation}
is called a \emph{first integral}. Equation \eqref{eq:def:CL} is
said to be a \emph{conservation law}.
\end{definition}
Equation \eqref{eq:preLC} asserts that $C = H T -
\boldsymbol{\psi}\emph{\textsuperscript{T}} \cdot \mathbf{X}$ is a
first integral:
\begin{theorem}[Noether's theorem]
\label{thm:TNoether}
If (\ref{eq:transf2}) is a \emph{symmetry} of problem
(\ref{eq:funcionalCO})-(\ref{eq:sistCont}), then
\begin{equation}
\label{eq:leicons} \boldsymbol{\psi}(t)\emph{\textsuperscript{T}}
\cdot \mathbf{X}
(t,\mathbf{x}(t),\psi_0,\boldsymbol{\psi}(t))-H(t,\mathbf{x}(t),
\mathbf{u}(t),\psi_0,\boldsymbol{\psi}(t))T(t,\mathbf{x}(t),\psi_0,\boldsymbol{\psi}(t))
= \textrm{const}
\end{equation}
is a conservation law.
\end{theorem}
From expression \eqref{eq:leicons} we see that Noetherian
conservation laws, associated with a certain optimal control
problem, that is, with a certain Hamiltonian
$H(t,\mathbf{x},\mathbf{u},\psi_0,\boldsymbol{\psi})$, only depend
on the generators $T$ and $\mathbf{X}$ of a symmetry
$\left(T,\mathbf{X},\mathbf{U},\boldsymbol{\Psi}\right)$
\eqref{eq:transf2}.


\section{Computation of Conservation Laws}
\label{sec:CCL}

The conservation laws we are looking for are obtained by
substitution of the components $T$ and $\mathbf{X}$ of a symmetry
of the problem in (\ref{eq:leicons}). In
section~\ref{sec:ProcMaple} we introduce the \maple procedure
\emph{Noether} to do these calculations for us. The input to this
procedure is: the Lagrangian $L$ and the velocity vector
$\boldsymbol{\varphi}$, that define the optimal control problem
\eqref{eq:funcionalCO}-\eqref{eq:sistCont} and the respective
Hamiltonian $H$; and a symmetry, or a family of symmetries,
obtained by means of our procedure \emph{Symmetry} (see the
procedure \emph{Symmetry} in \S\ref{sec:ProcMaple}). The output of
\emph{Noether} is the corresponding conservation law
\eqref{eq:leicons}. The non-trivial part of the computation lies
in the determination of the symmetries of the problem (implemented
in the \maple procedure \emph{Symmetry}). Our algorithm for
determining the infinitesimal generators is based on the necessary
and sufficient condition of invariance \eqref{eq:detGerad}. The
key to do the calculations consists in observing that when we
substitute the Hamiltonian $H$ and its partial derivatives in the
invariance identity \eqref{eq:detGerad}, then the condition
becomes a polynomial in $\dot{\mathbf{x}}$ and
$\dot{\boldsymbol{\psi}}$, and one can set the coefficients equal
to zero. Let us see how it works in detail.

Substituting $H$ and its partial derivatives into
\eqref{eq:detGerad}, and expanding the total derivatives
\[
\frac{\textrm{d}T}{\textrm{d}t} = \frac{\partial T}{\partial t} +
\frac{\partial T}{\partial \mathbf{x}}\cdot\dot{\mathbf{x}} +
\frac{\partial T}{\partial
\boldsymbol{\psi}}\cdot\dot{\boldsymbol{\psi}}\textrm{,}\quad
\frac{\textrm{d}\mathbf{X}}{\textrm{d}t} = \frac{\partial
\mathbf{X}}{\partial t} + \frac{\partial \mathbf{X}} {\partial
\mathbf{x}}\cdot\dot{\mathbf{x}}+ \frac{\partial
\mathbf{X}}{\partial \boldsymbol{\psi}}\cdot
\dot{\boldsymbol{\psi}},
\]
one can write equation (\ref{eq:detGerad}) as a polynomial
\begin{equation}
\label{eq:poly}
A(t,\mathbf{x},\mathbf{u},\psi_0,\boldsymbol{\psi}) +
B(t,\mathbf{x},\mathbf{u},\psi_0,\boldsymbol{\psi}) \cdot
\dot{\mathbf{x}} +
C(t,\mathbf{x},\mathbf{u},\psi_0,\boldsymbol{\psi})  \cdot
\dot{\boldsymbol{\psi}} = 0
\end{equation}
in the $2n$ derivatives $\dot{\mathbf{x}}$ and
$\dot{\boldsymbol{\psi}}$:
\begin{multline}
\label{eq:detGerad2} \left(\frac{\partial H}{\partial t}T
+\frac{\partial H}{\partial \mathbf{x}}\cdot \mathbf{X}
+\frac{\partial H}{\partial \mathbf{u}}\cdot \mathbf{U}
+\frac{\partial H}{\partial \boldsymbol{\psi}}\cdot
\boldsymbol{\Psi} +H \frac{\partial T}{\partial t}
-\boldsymbol{\psi}\textsuperscript{T} \cdot \frac{\partial
\mathbf{X}}{\partial t}\right) \\
+ \left( -\boldsymbol{\Psi}\textsuperscript{T} +H \frac{\partial
T}{\partial\mathbf{x}} -\boldsymbol{\psi}\textsuperscript{T} \cdot
\frac{\partial \mathbf{X}} {\partial \mathbf{x}} \right)\cdot
\dot{\mathbf{x}}+ \left( H \frac{\partial T}{\partial
\boldsymbol{\psi}} -\boldsymbol{\psi}\textsuperscript{T} \cdot
\frac{\partial \mathbf{X}}{\partial \boldsymbol{\psi}} \right)
\cdot \dot{\boldsymbol{\psi}}=0 \, .
\end{multline}
The terms in \eqref{eq:detGerad2} which involve derivatives with
respect to vectors are expanded in row-vectors or in matrices,
depending, respectively, if the function is a scalar or a
vectorial one. For example,
\begin{eqnarray}
\frac{\partial T}{\partial \mathbf{x}}
&=&\left[ \frac{\partial T}{\partial x_1}
\ \frac{\partial T}{\partial x_2}
\ \cdots \ \frac{\partial T}{\partial x_n}\right],
\nonumber\\
\frac{\partial \mathbf{X}}{\partial \boldsymbol{\psi}}&=&
\left[\frac{\partial \mathbf{X}}{\partial \psi_1}\ \frac{\partial
\mathbf{X}}{\partial \psi_2}\ \cdots\  \frac{\partial
\mathbf{X}}{\partial \psi_n}\right] = \left[
\begin{array}{cccc}
\smallskip
\frac{\partial X_1}{\partial \psi_1} & \frac{\partial X_1}{\partial \psi_2}
& \cdots & \frac{\partial X_1}{\partial \psi_n}\\
\frac{\partial X_2}{\partial \psi_1} & \frac{\partial X_2}{\partial \psi_2}
& \cdots & \frac{\partial X_2}{\partial \psi_n}\\
\vdots&\vdots&\ddots&\vdots\\
\frac{\partial X_n}{\partial \psi_1} & \frac{\partial X_n}{\partial \psi_2}
& \cdots & \frac{\partial X_n}{\partial \psi_n}
\end{array}
\right].\nonumber
\end{eqnarray}
Given an optimal control problem, defined by a Lagrangian $L$ and
a control system (\ref{eq:sistCont}), we determine the
infinitesimal generators $T$, $\mathbf{X}$, $\mathbf{U}$ and
$\boldsymbol{\Psi}$, which define a symmetry for the problem, by
the following method. Equation (\ref{eq:detGerad2}) is a
differential equation in the $2n+m+1$ unknown functions $T$,
$X_1$, \dots, $X_n$, $U_1$, \dots, $U_m$, $\varPsi_1$, \dots, and
$\varPsi_n$. This equation must hold for all $\dot{x}_1$, \dots,
$\dot{x}_n$, $\dot{\psi}_1$, \dots, $\dot{\psi}_n$, and therefore
the coefficients $A$, $B$, and $C$ of polynomial \eqref{eq:poly}
must vanish, that is,
\begin{equation}
\label{eq:detGerad3}
\left\{
\begin{array}{l}
\displaystyle \medskip
\frac{\partial H}{\partial t}T
+\frac{\partial H}{\partial \mathbf{x}}\cdot \mathbf{X}
+\frac{\partial H}{\partial \mathbf{u}}\cdot \mathbf{U}
+\frac{\partial H}{\partial \boldsymbol{\psi}}\cdot \boldsymbol{\Psi}
+H \frac{\partial T}{\partial t}
-\boldsymbol{\psi}\textsuperscript{T} \cdot
\frac{\partial \mathbf{X}}{\partial t}=0 \, , \\
\displaystyle \medskip
-\boldsymbol{\Psi}\textsuperscript{T}
+H \frac{\partial T}{\partial\mathbf{x}}
-\boldsymbol{\psi}\textsuperscript{T} \cdot \frac{\partial \mathbf{X}}
{\partial \mathbf{x}}
=\mathbf{0} \, , \\
\displaystyle H \frac{\partial T}{\partial \boldsymbol{\psi}}
-\boldsymbol{\psi}\textsuperscript{T} \cdot \frac{\partial
\mathbf{X}}{\partial \boldsymbol{\psi}} =\mathbf{0} \, .
\end{array}
\right.
\end{equation}
System of equations \eqref{eq:detGerad3}, obtained from
(\ref{eq:detGerad2}), is a system of $2n+1$ partial differential
equations with $2n+m+1$ unknown functions (so, in general, there
exists not a unique symmetry but a whole family of symmetries --
see examples in \S\ref{sec:EI} and \S\ref{sec:CV}). Although a
system of PDEs, its solution is possible because the system is
linear with respect to the unknown functions and their
derivatives. However, when dealing with optimal control problems
with several state and control variables, the number of
calculations is big enough, and the help of the computer is more
than welcome. Our \maple procedure \emph{Symmetry}, in
\S\ref{sec:ProcMaple}, does the job for us. The procedure
receives, as input, the Lagrangian and the expressions which
define the control system; and gives, as output, a family of
symmetries
$\left(T,\mathbf{X},\mathbf{U},\boldsymbol{\Psi}\right)$. Since
system \eqref{eq:detGerad3} is homogeneous, we always have, as
trivial solution,
$\left(T,\mathbf{X},\mathbf{U},\boldsymbol{\Psi}\right) =
\mathbf{0}$. This does not give any additional information (for
the trivial solution, Noether's theorem is the truism ``zero is a
constant''). When the output of \emph{Symmetry} coincides with the
trivial solution, that means, roughly speaking, that the optimal
control problem does not admit a symmetry (more precisely -- see
\S\ref{sec:ProcMaple} -- it means that our procedure \maple was
not able to find symmetries for the problem).

Summarizing: given an optimal control problem
\eqref{eq:funcionalCO}-\eqref{eq:sistCont} we compute conservation
laws, in an automatic way, through two steps: (i) with our
procedure \emph{Symmetry} we obtain all the possible invariance
symmetries of the problem; (ii) using the obtained symmetries as
input to procedure \emph{Noether}, based on
Theorem~\ref{thm:TNoether}, we obtain the correspondent
conservation laws.

In the next two sections we give several examples illustrating the
whole process.


\section{Illustrative Examples}
\label{sec:EI}

In order to show the functionality and the use of the routines developed,
we apply our \maple package to several concrete optimal
control problems found in the literature. The obtained results
show the correctness and usefulness of the \maple code.
All the computational processing was carried
out with \maple version 9.5 on a
1.4 GHz Pentium Centrino with 512MB RAM.
The computing time of procedure \emph{Symmetry}
is indicated for each example in the format min'sec''.
All the other \maple commands run instantaneously.


\begin{example}
\label{ex:Ex1}
(0'02'') Let us begin with the minimization of the
functional $\int_a^b L(u(t))\mathrm{d}t$ subject to the control
system $\dot x(t)=\varphi (u(t))x(t)$. This is a very simple
problem, with one state variable ($n = 1$) and one control
variable ($m = 1$). With \maple definitions
\small
\begin{verbatim}
> l:=L(u); Phi:=phi(u)*x;
\end{verbatim}
\begin{eqnarray*}
l&:=&L \left( u \right) \\
\Phi&:=&\varphi \left( u \right) x
\end{eqnarray*}
\normalsize
our procedure \emph{Symmetry} determines
the infinitesimal invariance generators of the optimal control
problem under consideration:
\small
\begin{verbatim}
> Symmetry(l,Phi,t,x,u);
\end{verbatim}
\[
\left\{ T=C_{2},\, X=C_{1}x,\,U=0,\,\Psi=-C_{1}\,\psi \right\}
\]
\normalsize
The family of Conservation Laws associated with the generators just
obtained, is easily obtained through our procedure
\emph{Noether} (the sign of percentage -- \% -- is an operator used in \maple
to represent the result of the previous command):
\small
\begin{verbatim}
> Noether(l,Phi,t,x,u, %);
\end{verbatim}
\[
C_{1}x(t) \psi(t) - \left( \psi_{0}L
 \left( u(t)  \right) +\psi(t) \varphi
 \left( u(t)  \right) x(t)  \right) C_{{2}
}={\it const}
\]
\normalsize
The obtained conservation law depends on two parameters.
Since the problem is autonomous,
the fact that the Hamiltonian $H = \psi_{0}L
 \left( u(t)  \right) +\psi(t) \varphi
 \left( u(t)  \right) x(t)$ is constant
along the extremals is a trivial consequence of
the property \eqref{eq:prop:dHdt}.
With the substitutions
\small
\begin{verbatim}
> subs(C[1]=1,C[2]=0, %);
\end{verbatim}
\[
x(t) \psi(t) ={\it const}
\]
\normalsize
one gets the conservation law obtained in \cite[Example~4]{delfimIO}.
\end{example}


\begin{example}
\label{ex:Ex2}
(1'13'') Let us consider now the following problem:
\begin{eqnarray*}
\int_a^b \left(u_1(t)^2+u_2(t)^2\right)\mathrm{d}t \longrightarrow \min,\\
\left\{\begin{array}{l}
\dot x_1(t)=u_1(t) \cos x_3(t),\\
\dot x_2(t)=u_1(t) \sin x_3(t),\\
\dot x_3(t)=u_2(t),
\end{array}\right.
\end{eqnarray*}
where the control system serves as model for the kinematics of a
car \cite[Example~18, p.~750]{Rouchon01}. In this case the optimal
control problem has three state variables ($n=3$) and two controls
($m=2$). The conservation law for this example, and the next ones,
is obtained by the same process followed in Example~\ref{ex:Ex1}.
\small
\begin{verbatim}
> L:=u[1]^2+u[2]^2; phi:=[u[1]*cos(x[3]),u[1]*sin(x[3]),u[2]];
\end{verbatim}
\begin{eqnarray*}
L&:=&{u_{1}}^{2}+{u_{2}}^{2}\\
\varphi&:=&[u_{1}\cos \left( x_{3} \right),
u_{1}\sin \left( x_{3} \right) ,u_{2}]
\end{eqnarray*}
\normalsize
\small
\begin{verbatim}
> Symmetry(L, phi, t, [x[1],x[2],x[3]], [u[1],u[2]]);
\end{verbatim}
\begin{multline*}
\left\{ T=C_{2},\, X_{1}=-C_{1}x_{2}+C_{3},\,X_{2}=C_{1}x_{1}+C_{4},\,X_{3}=C_{1},\right.\\
\left.U_{1}=0,\,U_{2}=0,\,
\Psi_{1}=-C_{1}\psi_{2},\,\Psi_{2}=C_{1}\psi_{1},\,\Psi_{3}=0\right\}
\end{multline*}
\normalsize
\small
\begin{verbatim}
> Noether(L, phi, t, [x[1],x[2],x[3]], [u[1],u[2]], %);
\end{verbatim}
\begin{multline*}
{\left( -C_{1}x_{2}(t) +C_{3} \right)
\psi_{1}(t) + \left( C_{1}x_{1}(t)
+C_{4} \right) \psi_{2}(t)
\mbox{}+C_{1}\psi_{3}(t)}\\
-\biggl( \psi_{0} \left(  \left( u_{1}(t)  \right) ^{2}
+ \left( u_{2}(t)  \right) ^{2} \right)
\hbox{}+u_{1}(t) \cos
\left( x_{3}(t)  \right) \psi_{1}(t)
\mbox{}+u_{1}(t) \sin \left( x_{3}(t)\right)
\psi_{2}(t) +u_{2}(t) \psi_{3}
(t)\biggr) C_{2}={\it const}
\end{multline*}
\normalsize
Choosing $C_1=1$ and $C_2=C_3=C_4=0$ we obtain,
from Theorem~\ref{thm:TNoether}, the conservation law
\small
\begin{verbatim}
> subs(C[1]=1,C[2]=0,C[3]=0,C[4]=0, %);
\end{verbatim}
\[
-x_{2}(t) \psi_{1}(t)
+x_{1}(t) \psi_{2}(t)
+\psi_{3}(t) ={\it const}
\]
\normalsize
which corresponds to the symmetry group of
planar (orientation-preserving) isometries
given in \cite[Example~18, p.~750]{Rouchon01}.
\end{example}


\begin{example}
\label{ex:Ex3}
(0'01'')
Let us return to a scalar problem ($n=m=1$):
\begin{equation*}
\int_a^b \mathrm{e}^{t x(t)} u(t)\mathrm{d}t \longrightarrow \min, \quad
\dot x(t)=t x(t) u(t)^2 \, .
\end{equation*}
\small
\begin{verbatim}
> L:=exp(t*x)*u; phi:=t*x*u^2;
\end{verbatim}
\begin{eqnarray*}
L&:=&{e^{tx}}u\\
\varphi&:=&tx{u}^{2}
\end{eqnarray*}
\normalsize
\small
\begin{verbatim}
> Symmetry(L, phi, t, x, u);
\end{verbatim}
\[
 \left\{ T=-tC_{1},\, X=C_{1}x,\,U=C_{1}u,\,\Psi=-\psi\,C_{1}\right\}
\]
\normalsize
\small
\begin{verbatim}
> Noether(L, phi, t, x, u, %);
\end{verbatim}
\[
C_{1}x(t) \psi(t)
+ \left( \psi_{0}{e^{tx(t) }}u(t)
+\psi(t) tx(t)  \left( u(t)
\right) ^{2} \right) tC_{1}={\it const}
\]
\normalsize
By choosing $C_1=1$
\small
\begin{verbatim}
> expand(subs(C[1]=1, %));
\end{verbatim}
\[
x(t) \psi(t) +t\psi_{0}{e^{tx(t)}}u(t)
+\psi(t) {t}^{2}x(t)\left( u(t)  \right) ^{2}={\it const}
\]
\normalsize
one obtains the conservation law of \cite[Example~1]{3ncnw}.
\end{example}


\begin{example}
\label{ex:Ex5}
(6'41'')
We now consider an optimal control problem with four state variables ($n=4$)
and two controls ($m=2$):
\begin{eqnarray*}
&&\int_a^b \left(u_1(t)^2+u_2(t)^2\right)\mathrm{d}t \longrightarrow \min,\\
&&\left\{\begin{array}{l}
\smallskip
\dot x_1(t)=x_3(t),\\
\smallskip
\dot x_2(t)=x_4(t),\\
\smallskip
\dot x_3(t)=-x_1(t)\left(x_1(t)^2+x_2(t)^2\right)+u_1(t),\\
\dot x_4(t)=-x_2(t)\left(x_1(t)^2+x_2(t)^2\right)+u_2(t),
\end{array}\right.
\end{eqnarray*}
\small
\begin{verbatim}
> L:=u[1]^2+u[2]^2; phi:=[x[3],x[4],-x[1]*(x[1]^2+x[2]^2)+u[1],
     -x[2]*(x[1]^2+x[2]^2)+u[2]];
\end{verbatim}
\begin{eqnarray*}
L&:=&{u_{1}}^{2}+{u_{2}}^{2}\\
\varphi&:=&\left[x_{3},x_{4},-x_{1} \left( {x_{1}}^{2}+{x_{2}}^{2}
\right) +u_{1},-x_{2} \left( {x_{1}}^{2}+{x_{2}}^{2} \right)
+u_{2}\right]
\end{eqnarray*}
\normalsize
\small
\begin{verbatim}
> Symmetry(L, phi, t, [x[1],x[2],x[3],x[4]], [u[1],u[2]]);
\end{verbatim}
\begin{multline*}
\Bigl\{T=C_{3},\,
X_{1}=C_{1}x_{2},\,X_{2}=-C_{1}x_{1},\,X_{3}=C_{1}x_{4},X_{4}=-C_{1}\mbox{}x_{3},\,\\
U_{1}=-C_{2}u_{2}\mbox{}-\frac{1}{2}\,{\frac { \left( C_{2}+C_{1} \right) \psi_{4}}{\psi_{0}}},\,
U_{2}=C_{2}u_{1}+\frac{1}{2}\,{\frac { \left( C_{2}+C_{1}\right) \psi_{3}}{\psi_{0}}},\,\\
\Psi_{1}=C_{1}\psi_{2},\,\Psi_{2}=-C_{1}\psi_{1},\,
\Psi_{3}=C_{1}\mbox{}\psi_{4},\,\Psi_{4}=-C_{1}\psi_{3}\Bigr\}
\end{multline*}
\normalsize
\small
\begin{verbatim}
> Noether(L, phi, t, [x[1],x[2],x[3],x[4]], [u[1],u[2]], %);
\end{verbatim}
\begin{eqnarray*}
C_{1}x_{2}(t) \psi_{1}(t) -C_{1}x_{1}
(t) \psi_{2}(t)
\mbox{}+C_{1}x_{4}(t) \psi_{3}(t)
-C_{1}x_{3}(t) \psi_{4}(t)
- \biggl( \psi_{0} \left(  \left( u_{1}(t)  \right) ^{2}
+ \left( u_{2}(t)  \right) ^{2} \right)\\ \hbox{}
+x_{3}(t) \psi_{1}(t)
+x_{4}(t) \psi_{2}(t)
\mbox{}+ \left( -x_{1}(t)  \left(  \left( x_{1}
(t)  \right) ^{2}+ \left( x_{2}(t)  \right) ^{2}
\right) +u_{1}(t)  \right) \psi_{3}(t)\\ \hbox{}
+ \left( -x_{2}(t)  \left(  \left( x_{1}(t)
\right) ^{2}+ \left( x_{2}(t)  \right) ^{2} \right) +u_{2}
(t)  \right) \psi_{4}(t)  \biggr) C_{3}
\mbox{}={\it const}
\end{eqnarray*}
\normalsize
The substitutions
\small
\begin{verbatim}
> subs(C[1]=-1,C[3]=0, %);
\end{verbatim}
\[
-x_{2}(t) \psi_{1}(t) +x_{1}\left( t\right)
\psi_{2}(t) -x_{4} \left(t\right)\psi_{3}\left(t\right)
\mbox{}+x_{3}(t) \psi_{4}(t) ={\it const}
\]
\normalsize
conduce us to the conservation law provided in \cite[Example 5.2]{Torres04}.
\end{example}


\begin{example}
\label{ex:Ex6}
(6'42'')
Another problem with $n=4$ and $m=2$:
\begin{eqnarray*}
&&\int_a^b \left(u_1(t)^2+u_2(t)^2\right)\mathrm{d}t \longrightarrow \min,\\
&&\left\{\begin{array}{l}
\smallskip
\dot x_1(t)=u_1(t)\left(1+x_2(t)\right),\\
\smallskip
\dot x_2(t)=u_1(t)x_3(t),\\
\smallskip
\dot x_3(t)=u_2(t),\\
\dot x_4(t)=u_1(t)x_3(t)^2.
\end{array}\right.
\end{eqnarray*}
\small
\begin{verbatim}
> L:=u[1]^2+u[2]^2; phi:=[u[1]*(1+x[2]),u[1]*x[3],u[2],u[1]*x[3]^2];
\end{verbatim}
\begin{eqnarray*}
L&:=&{u_{1}}^{2}+{u_{2}}^{2}\\
\varphi &:=&[u_{1} \left( 1+x_{2} \right) ,u_{1}x_{3},u_{2},
u_{1}{x_{3}}^{2}]
\end{eqnarray*}
\normalsize
\small
\begin{verbatim}
> Symmetry(L, phi, t, [x[1],x[2],x[3],x[4]], [u[1],u[2]]);
\end{verbatim}
\begin{multline*}
\Bigl\{T=\frac{2}{3}\,C_{1}t+C_{2},\,
X_{1}=C_{1}x_{1}+C_{4},\, X_{2}=\frac{2}{3}\,C_{1}+\frac{2}{3}\,C_{1}x_{2},\,
X_{3}=\frac{1}{3}\,C_{1}x_{3},\,  X_{4}=C_{1}x_{4}+C_{3},\,\\
U_{1}=-\frac{1}{3}\,u_{1}C_{1},\, U_{2}=-\frac{1}{3}\,C_{1}u_{2},\,
\Psi_{1}=-C_{1}\psi_{1},\, \Psi_{2}=-\frac{2}{3}\, C_{1}\psi_{2},\,
\Psi_{3}=-\frac{1}{3}\,C_{1}\psi_{3},\, \Psi_{4}=-C_{1}\psi_{4}\Bigr\}
\end{multline*}
\normalsize
\small
\begin{verbatim}
> Noether(L, phi, t, [x[1],x[2],x[3],x[4]], [u[1],u[2]], %);
\end{verbatim}
\begin{eqnarray*}
\left( C_{1}x_{1}(t) +C_{4} \right) \psi_{1}
(t) + \left( \frac{2}{3}\,C_{1}+\frac{2}{3}\,C_{1}x_{2}(t)
\right) \psi_{2}(t)
\mbox{}+\frac{1}{3}\,C_{1}x_{3}(t) \psi_{3}(t)
\mbox{}+ \left( C_{1}x_{4}(t) +C_{3} \right) \psi_{4}
(t) \\ \hbox{}
- \Biggl( \psi_{0} \left(  \left( u_{1}(t)  \right) ^{2}
+ \left( u_{2}(t)  \right) ^{2} \right) +\psi_{1}
(t) u_{1}(t)  \left( 1+x_{2}
(t)  \right)+\psi_{2}(t) u_{1}(t)
x_{3}(t) \\ \hbox{}
+\psi_{3}(t) u_{2}(t) +\psi_{4}
(t) u_{1}(t)  \left( x_{3}(t)
\right) ^{2} \Biggr)  \left( \frac{2}{3}\,C_{1}t+C_{2} \right)
\mbox{}={\it const}
\end{eqnarray*}
\normalsize
With the substitutions
\small
\begin{verbatim}
> subs(C[1]=3,C[2]=0,C[3]=0,C[4]=0, %);
\end{verbatim}
\begin{eqnarray*}
3\,x_{1}(t) \psi_{1}(t) + \left( 2+2\,x_{2}
(t)  \right) \psi_{2}(t) +x_{3}(t)
\psi_{3}(t)+3\,x_{4}(t) \psi_{4}
(t)-2\, \Bigl( \psi_{0} \left(  \left( u_{1}(t)
\right) ^{2}+ \left( u_{2}(t)  \right) ^{2} \right) \\ \hbox{}
+\psi_{1}(t) u_{1}(t)  \left( 1+x_{2}
(t)  \right)+\psi_{2}(t) u_{1}(t)
x_{3}(t) +\psi_{3}(t) u_{2}(t)
+\psi_{4}(t) u_{1}(t)  \left( x_{3}
(t)  \right) ^{2} \Bigr) t={\it const}
\end{eqnarray*}
\normalsize
we have the conservation law obtained in \cite[Example 5.3]{Torres04}.
\end{example}


\begin{example}
\label{ex:Ex7}
(0'04'')
Let us consider
\begin{eqnarray*}
&&\int_a^b u(t)^2\mathrm{d}t \longrightarrow \min,\\
&&\left\{\begin{array}{l}
\dot x(t)=1+y(t)^2,\\
\dot y(t)=u(t).
\end{array}\right.
\end{eqnarray*}
\small
\begin{verbatim}
> L:=u^2; phi:=[1+y^2,u];
\end{verbatim}
\begin{eqnarray*}
L&:=&{u}^{2}\\
\varphi &:=&[1+{y}^{2},u]
\end{eqnarray*}
\normalsize
\small
\begin{verbatim}
> Symmetry(L, phi, t, [x,y], u);
\end{verbatim}
\begin{eqnarray*}
\left\{T=\frac{1}{2}\,C_{1}t+C_{2},\,
X_{1}=-\frac{1}{2}\,C_{1}t+C_{1}x+C_{3}\mbox{},\,
X_{2}=\frac{1}{4}\,C_{1}y,\,\right.\\
\left. U=-\frac{1}{4} C_{1}u,\,
\Psi_{1}=-C_{1}\psi_{1},\,
\Psi_{2}=-\frac{1}{4}\,C_{1}\psi_{2}\right\}
\end{eqnarray*}
\normalsize
\small
\begin{verbatim}
> Noether(L, phi, t, [x,y], u, %);
\end{verbatim}
\begin{multline*}
\left( -\frac{1}{2}\,C_{1}t+C_{1}x(t) +C_{3} \right) \psi_{1}
(t) +\frac{1}{4}\,C_{1}y(t) \psi_{2}(t)\\
- \Bigl( \psi_{0} \left( u(t)  \right) ^{2}+\psi_{1}
(t)  \left( 1+ \left( y(t)  \right) ^{2} \right)
+\psi_{2}(t) u(t)  \Bigr)  \left( \frac{1}{2}\,C_{1}t
+C_{2} \right)={\it const}
\end{multline*}
\normalsize
From substitutions
\small
\begin{verbatim}
> subs(C[1]=-4,C[2]=0,C[3]=0,%);
\end{verbatim}
\begin{eqnarray*}
\left( 2\,t-4\,x(t)  \right) \psi_{1}(t)
-y(t) \psi_{2}(t) +2\, \left( \psi_{0}
\left( u(t)  \right) ^{2}+\psi_{1}(t)
\left( 1+ \left( y(t)  \right) ^{2} \right) +\psi_{2}
(t) u(t)  \right) t ={\it const}
\end{eqnarray*}
\normalsize
we obtain the conservation law in \cite[Example 6.2]{Torres04}.
\end{example}


\begin{example}
\label{ex:Ex8}
(2'44'')
We consider now a minimum time problem ($T \rightarrow \min
\Leftrightarrow \int_0^T 1\mathrm{d}t \rightarrow \min$)
with the following control system:
\[
\left\{\begin{array}{l}
\dot x_1(t)=1+x_2(t),\\
\dot x_2(t)=x_3(t),\\
\dot x_3(t)=u(t),\\
\dot x_4(t)=x_3(t)^2-x_2(t)^2 .
\end{array}\right.
\]
In \maple we have:
\small
\begin{verbatim}
> L:=1; phi:=[1+x[2],x[3],u,x[3]^2-x[2]^2];
\end{verbatim}
\begin{eqnarray*}
L&:=&1\\
\varphi &:=&[1+x_{2},x_{3},u,{x_{3}}^{2}-{x_{2}}^{2}]
\end{eqnarray*}
\normalsize
\small
\begin{verbatim}
> Symmetry(L, phi, t, [x[1],x[2],x[3],x[4]], u);
\end{verbatim}
\begin{multline*}
\left\{T=C_{5},\,
X_{1}= \left( -\frac{1}{2}\,C_{2}-\frac{1}{2}\,C_{1} \right) t +\frac{1}{2}\,C_{2}x_{1}+C_{4},\,
X_{2}=-\frac{1}{2}\,C_{1}+\frac{1}{2}\,C_{2}x_{2},\,\right.\\
X_{3}=\frac{1}{2}\,C_{2}x_{3},\, X_{4}=-C_{1}t+C_{1}x_{1}+C_{2}x_{4}+C_{3},\,\\
\left. U=\frac{1}{2}\,uC_{2},\,
\Psi_{1}=-\frac{1}{2}\,C_{2}\psi_{1}-C_{1}\psi_{4}, \, \Psi_{2}=-\frac{1}{2}\,C_{2}\psi_{2}, \,
\Psi_{3}=-\frac{1}{2}\,C_{2}\psi_{3}, \, \Psi_{4}=-C_{2}\psi_{4} \right\}
\end{multline*}
\normalsize
\small
\begin{verbatim}
> Noether(L, phi, t, [x[1],x[2],x[3],x[4]], u, %);
\end{verbatim}
\begin{multline*}
\left(  \left( -\frac{1}{2}\,C_{2}-\frac{1}{2}\,C_{1} \right) t
+\frac{1}{2}\,C_{2}x_{1}
(t)+C_{4} \right) \psi_{1}(t) + \left( -\frac{1}{2}\,
C_{1}+\frac{1}{2}\,C_{2}x_{2}(t)  \right) \psi_{2}(t)\\
+\frac{1}{2}\,C_{2}x_{3}(t) \psi_{3}
(t) + \left( -C_{1}t+C_{1}x_{1}(t) +C_{2}x_{4}
(t) +C_{3} \right) \psi_{4}(t)\\
-\left( \psi_{0}+ \left( 1+x_{2}(t)  \right) \psi_{1}(t)
+x_{3}(t) \psi_{2}(t)+\psi_{3}(t) u(t)
+ \left(  \left( x_{3}(t)  \right) ^{2}
- \left( x_{2}(t)\right)^{2} \right) \psi_{4}(t)\right) C_{5}={\it const}
\end{multline*}
\normalsize
With the appropriate values for the constants,
\small
\begin{verbatim}
> subs(C[1]=0,C[2]=2,C[3]=0,C[4]=0,C[5]=0,%);
\end{verbatim}
\[
 \left( -t+x_{1}(t)  \right) \psi_{1}(t)
 +x_{2}(t) \psi_{2}(t) +x_{3}
(t) \psi_{3}(t)
\mbox{}+2\,x_{4}(t) \psi_{4}(t) ={\it const}
\]
\normalsize
we obtain the conservation law derived in \cite[Example 6.3]{Torres04}.
\end{example}


\begin{example}
\label{ex:Ex9}
(0'25'')
Follows another problem of minimum time, with
control system given by
\[
\left\{\begin{array}{l}
\dot x(t)=1+y(t)^2-z(t)^2\, ,\\
\dot y(t)=z(t)\, ,\\
\dot z(t)=u(t)\, .
\end{array}\right.
\]
\small
\begin{verbatim}
> L:=1; phi:=[1+y^2-z^2,z,u];
\end{verbatim}
\begin{eqnarray*}
L&:=&1\\
\varphi &:=&[1+{y}^{2}-{z}^{2},z,u]
\end{eqnarray*}
\normalsize
\small
\begin{verbatim}
> Symmetry(L, phi, t, [x,y,z], u);
\end{verbatim}
\begin{multline*}
\left\{ T=C_{2},\,
X_{1}=-C_{1}t+C_{1}x+C_{3},\, X_{2}=\frac{1}{2}\,C_{1}y,\, X_{3}=\frac{1}{2}\,C_{1}z,\, \right.\\
\left. U=\frac{1}{2}\,C_{1}u,\,
\Psi_{1}=-C_{1}\psi_{1},\, \Psi_{2}=-\frac{1}{2}\,C_{1}\psi_{2},\,
\Psi_{3}= \mbox{}-\frac{1}{2}\,C_{1}\psi_{3}\right\}
\end{multline*}
\normalsize
\small
\begin{verbatim}
> Noether(L, phi, t, [x,y,z], u, %);
\end{verbatim}
\begin{eqnarray*}
\left( -C_{1}t+C_{1}x(t) +C_{3} \right) \psi_{1}
(t)+\frac{1}{2}\,C_{1}y(t) \psi_{2}(t)
+\frac{1}{2}\,C_{1}z(t) \psi_{3}(t)\\ \hbox{}
- \left( \psi_{0}+\psi_{1}(t)  \left( 1+
\left( y(t)  \right) ^{2}- \left( z(t)
\right) ^{2} \right)+\psi_{2}(t) z(t) +\psi_{3}
(t) u(t)\right) C_{2}={\it const}
\end{eqnarray*}
\normalsize
Substitutions
\small
\begin{verbatim}
> subs(C[1]=2,C[2]=0,C[3]=0, %);
\end{verbatim}
\[
 \left( -2\,t+2\,x(t)  \right) \psi_{1}(t)
 +y(t) \psi_{2}(t) +z(t) \psi_{3}
(t) ={\it const}
\]
\normalsize
specify the conservation law in the one obtained in \cite[Example 6.4]{Torres04}.
\end{example}


We finish the section by applying our \maple package to three important
problems of geodesics in sub-Riemannian geometry.
The reader, interested in the study of symmetries of flat distributions
of sub-Riemannian geometry, is referred to \cite{Sachkov04}.


\begin{example}[Martinet -- $\mathbf{(2,2,3)}$ problem]
\label{ex:Ex4}
Given the problem ($n=3,m=2$)
\begin{eqnarray*}
&&\int_a^b \left(u_1(t)^2+u_2(t)^2\right)\mathrm{d}t \longrightarrow \min,\\
&&\left\{\begin{array}{l}
\smallskip
\dot x_1(t)=u_1(t),\\
\displaystyle \smallskip
\dot x_2(t)=\frac{u_2(t)}{1+\alpha x_1(t)},\\
\dot x_3(t)=x_2(t)^2 u_1(t),
\end{array}\right.\qquad \alpha \in \mathbb{R},
\end{eqnarray*}
we consider two distinct situations: $\alpha=0$ (Martinet problem
of sub-Riemannian geometry in the flat case -- see \cite{Bonnard98})
and $\alpha \ne 0$ (non-flat case).\\

\emph{Flat Problem ($\alpha=0$, 1'06''):}
\small
\begin{verbatim}
> L:=u[1]^2+u[2]^2; phi:=[u[1],u[2],x[2]^2*u[1]];
\end{verbatim}
\begin{eqnarray*}
L&:=&{u_{1}}^{2}+{u_{2}}^{2}\\
\varphi&:=&[u_{1},u_{2},{x_{2}}^{2}u_{1}]
\end{eqnarray*}
\normalsize
\small
\begin{verbatim}
> Symmetry(L, phi, t, [x[1],x[2],x[3]], [u[1],u[2]]);
\end{verbatim}
\begin{eqnarray*}
\left\{ T=\frac{2}{3}\,C_{1}t+C_{2},\,
X_{1}=\frac{1}{3}\,C_{1}x_{1}+C_{4},\,X_{2}=\frac{1}{3}\mbox{}\,C_{1}x_{2},\,
X_{3}=C_{1}x_{3}\mbox{}+C_{3},\, \right.\\
\left. U_{1}=-\frac{1}{3}\,u_{1}C_{1},\, U_{2}=-\frac{1}{3}\,C_{1}u_{2},\,
\Psi_{1}=-\frac{1}{3}\,C_{1}\psi_{1},\,
\Psi_{2}=-\frac{1}{3}\,C_{1}\psi_{2},\,
\Psi_{3}=-C_{1}\psi_{3}\right\}
\end{eqnarray*}
\normalsize
\small
\begin{verbatim}
> Noether(L, phi, t, [x[1],x[2],x[3]], [u[1],u[2]], %);
\end{verbatim}
\begin{multline*}
\left( \frac{1}{3}\,C_{1}x_{1}(t) +C_{4} \right) \psi_{1}(t)
+\frac{1}{3}\,C_{1}x_{2}(t) \psi_{2}(t)
+ \left( C_{1}x_{3}(t)
+C_{3} \right) \psi_{3}(t)\\
- \biggl( \psi_{0} \left(  \left(u_{1}(t)\right)^{2}
+\left( u_{2}(t)  \right) ^{2} \right)
 +\psi_{1}(t) u_{1}(t) +\psi_{2}(t) u_{2}(t)
 +\psi_{3}(t)  \left( x_{2}(t)\right)^{2}u_{1}(t)  \biggr)
\left( \frac{2}{3}\,C_{1}t+C_{2} \right) ={\it const}
\end{multline*}
\normalsize
With the substitutions
\small
\begin{verbatim}
> subs(C[1]=3,C[2]=0,C[3]=0,C[4]=0, %);
\end{verbatim}
\begin{multline*}
x_{1}(t) \psi_{1}(t) +x_{2}(t) \psi_{2}(t) +3\,x_{3}(t) \psi_{3}(t)\\
-2\, \Bigl( \psi_{0} \left(  \left( u_{1}(t)
\right) ^{2}+ \left( u_{2}(t)  \right) ^{2} \right)
+ \psi_{1}(t) u_{1}(t) +\psi_{2}(t) u_{2}(t)
 + \psi_{3}(t)  \left( x_{2}(t)
\right) ^{2}u_{1}(t)  \Bigr) t={\it const}
\end{multline*}
\normalsize
we get the conservation law first obtained in \cite[Example~2]{3ncnw}.\\

\emph{Non-Flat Problem ($\alpha \ne 0$, 1'14''):}
\small
\begin{verbatim}
> L:=u[1]^2+u[2]^2; phi:=[u[1],u[2]/(1+alpha*x[1]),x[2]^2*u[1]];
\end{verbatim}
\begin{eqnarray*}
L := {u_{1}}^{2}+{u_{2}}^{2}\\
\varphi:=\left[u_{1},{\frac {u_{2}}{1+\alpha\,
x_{1}}},{x_{2}}^{2}u_{1}\right]
\end{eqnarray*}
\normalsize
\small
\begin{verbatim}
> simplify(Symmetry(L, phi, t, [x[1],x[2],x[3]], [u[1],u[2]]));
\end{verbatim}
\begin{eqnarray*}
\bigl\{T=2\,C_{7}t+C_{11},\,
X_{1}=C_{7}(\alpha^{-1}+ x_{1}),\, X_{2}=0,\, X_{3}=C_{7}x_{3}+C_{10},\\
U_{1}=-C_{7}u_{1},\, U_{2}=-C_{7}u_{2},\,
\Psi_{1}=-C_{7}\psi_{1},\, \Psi_{2}=0,\, \Psi_{3}=-C_{7}\psi_{3}\bigr\}
\end{eqnarray*}
\normalsize
\small
\begin{verbatim}
> Noether(L, phi, t, [x[1],x[2],x[3]], [u[1],u[2]], %);
\end{verbatim}
\begin{multline*}
C_{7}(\alpha^{-1}+ x_{1}(t)) \psi_{1}(t) + \left( C_{7}x_{3}
(t) +C_{10} \right) \psi_{3}(t)\\
- \biggl( \psi_{0} \left(  \left( u_{1}(t)\right)^{2}
+ \left( u_{2}(t)  \right) ^{2} \right)
+u_{1}(t) \psi_{1}(t) +{\frac {u_{2}
(t) \psi_{2}(t) }{1+\alpha\,x_{1}
(t) }}+ \left( x_{2}(t)  \right) ^{2}u_{1}
(t) \psi_{3}(t)  \biggr)
\left( 2\,C_{7}t+C_{11} \right) ={\it const}
\end{multline*}
\normalsize
When $C_7=1$ and $C_{10} = C_{11}=0$,
\small
\begin{verbatim}
> subs(C[7]=1,C[10]=0,C[11]=0, %);
\end{verbatim}
\[
\begin{array}{l}
\smallskip
\left( {\alpha}^{-1}+x_{1}(t)  \right) \psi_{1}(t)
+x_{3}(t) \psi_{3}(t)\\
\displaystyle
\mbox{}-2\, \left( \psi_{0} \left(  \left( u_{1}(t)
\right) ^{2}+ \left( u_{2}(t)  \right) ^{2} \right)
+u_{1}(t) \psi_{1}(t) +{\frac {u_{2}
(t) \psi_{2}(t) }{1+\alpha\,
x_{1}(t) }}+ \left( x_{2}(t)  \right) ^{2}
u_{1}(t) \psi_{3}(t)  \right) t={\it const}
\end{array}
\]
\normalsize
we obtain the conservation law proved in \cite[Example~2]{delfimEJC}.
\end{example}


\begin{example}[Heisenberg -- $\mathbf{(2,3)}$ problem]
\label{ex:Ex10}
(1'04'')
The Heisenberg $(2,3)$ problem can be formulated as follows:
\begin{eqnarray*}
&&\frac{1}{2}\,\int_a^b \left({u_1(t)}^2+{u_2(t)}^{2}\right)\mathrm{d}t \longrightarrow \min,\\
&&\left\{\begin{array}{l}
\dot x_1(t)=u_1(t)\, ,\\
\dot x_2(t)=u_2(t)\, ,\\
\dot x_3(t)=u_2(t) x_1(t)\, .\\
\end{array}\right.
\end{eqnarray*}
The problem was proved to be completely integrable
using three independent conservation laws \cite{Taimanov97}.
Such conservation laws can now be easily obtained with our \maple functions.
\small
\begin{verbatim}
> L:=1/2*(u[1]^2+u[2]^2);
> phi:=[u[1], u[2], u[2]*x[1]];
\end{verbatim}
\begin{eqnarray*}
L&:=&\frac{1}{2}\,{u_{{1}}}^{2}+\frac{1}{2}\,{u_{{2}}}^{2}\\
\varphi &:=&[u_1, u_2, u_2 x_1]\\
\end{eqnarray*}
\begin{verbatim}
> Symmetry(L, phi, t, [x[1],x[2],x[3]], [u[1],u[2]]);
\end{verbatim}
\begin{multline*}
\bigl\{T=2\,C_2t+C_5,\,
X_1=C_1+C_2{\it x_1},\, X_2=C_2{\it x_2}+C_3,\, X_3=C_1{\it x_2}+2\,C_2{\it x_3}+C_4,\\
 U_1=-C_2{\it u_1},\, U_2=-C_2{\it u_2},\,
\Psi_1=-C_2\psi_1,\, \Psi_2=-C_2\psi_2-\psi_3C_1,\, \Psi_3=-2\,C_2\psi_3  \bigr\}
\end{multline*}
\normalsize
\small
\begin{verbatim}
> CL:=Noether(L, phi, t, [x[1],x[2],x[3]], [u[1],u[2]], %);
\end{verbatim}
\begin{multline*}
 CL:=\left( C_{{1}}+C_{{2}}x_{{1}}(t)  \right) \psi_{{1}}(t)
 + \left( C_{{2}}x_{{2}}(t) +C_{{3}}
 \right) \psi_{{2}}(t) + \left( C_{{1}}x_{{2}} \left( t
 \right) +2\,C_{{2}}x_{{3}}(t) +C_{{4}} \right) \psi_{{3}}(t) \\
- \biggl( \psi_{{0}} \left( \frac{1}{2}\, \left( u_{{1}}(t)
\right)^{2}+\frac{1}{2}\, \left( u_{{2}}(t)
 \right) ^{2} \right) +\psi_{{1}}(t) u_{{1}} \left( t
 \right) +\psi_{{2}}(t) u_{{2}}(t)
 +\psi_{{3}}(t) u_{{2}}(t) x_{{1}} \left( t
 \right)  \biggr)  \left( 2\,C_{{2}}t+C_{{5}} \right) ={\it const}
\end{multline*}
\normalsize
We now want to eliminate the controls from the previous family
of conservation laws. We begin to define the Hamiltonian:
\small
\begin{verbatim}
> H:=-L+Vector[row]([psi[1](t), psi[2](t), psi[3](t)]).Vector(phi);
\end{verbatim}
\[
H:=-\frac{1}{2}\,{u_1}^{2}-\frac{1}{2}\,{u_2}^{2}+u_1\psi_1(t)+u_2\psi_2(t)+u_2x_1\psi_3(t)
\]
\normalsize
The stationary condition \eqref{eq:stat:cond} permits to obtain
the pair of controls $(u_1(t),u_2(t))$.
\small
\begin{verbatim}
> solve({diff(H,u[1])=0, diff(H,u[2])=0}, {u[1],u[2]}):
> subs(x[1]=x[1](t), u[1]=u[1](t), u[2]=u[2](t), %);
\end{verbatim}
\[
 \left\{ u_1(t)=\psi_1(t),u_2(t)=\psi_2(t)+x_1(t)\psi_3(t)\right\}
\]
\normalsize
It is not difficult to show that the problem does not admit abnormal
extremals, so one can choose, without any loss of generality, $\psi_0 = -1$.
\small
\begin{verbatim}
> CL:=subs(psi[0]=-1, %, CL);
\end{verbatim}
\begin{multline*}
 CL:=\left( C_{{1}}+C_{{2}}x_{{1}}(t)  \right) \psi_{{1}}(t)
 + \left( C_{{2}}x_{{2}}(t) +C_{{3}}
 \right) \psi_{{2}}(t) + \left( C_{{1}}x_{{2}} \left( t
 \right) +2\,C_{{2}}x_{{3}}(t) +C_{{4}} \right) \psi_{{3}}(t)\\
- \biggl( \frac{1}{2}\, \left( \psi_{{1}}(t)
 \right) ^{2}-\frac{1}{2}\, \left( \psi_{{2}}(t) +x_{{1}}(t) \psi_{{3}}(t)\right)^{2}
 +\psi_{{2}}(t)  \left( \psi_{{2}}(t) +x_{{1}} \left(
t \right) \psi_{{3}}(t)  \right)\\
+\psi_{{3}} \left( t
 \right)  \left( \psi_{{2}}(t) +x_{{1}} \left( t
 \right) \psi_{{3}}(t)  \right) x_{{1}} \left( t
 \right)  \biggr)  \left( 2\,C_{{2}}t+C_{{5}} \right) ={\it const}
\end{multline*}
\normalsize
It is easy to extract from the family of conservation laws just obtained,
three independent conservation laws. We just need
to fix one constant to a non-zero value, and choose all the other constants to be zero:
\small
\begin{verbatim}
> subs(C[4]=1,seq(C[i]=0,i=1..5), CL);
> subs(C[1]=1,seq(C[i]=0,i=1..5), CL);
> simplify(subs(C[5]=-1,seq(C[i]=0,i=1..5), CL));
\end{verbatim}
\begin{eqnarray*}
&&\psi_{{3}}(t) ={\it const}\\
&&x_{{2}}(t) \psi_{{3}}(t) +\psi_{{1}}(t) = {\it const}\\
&&\frac{1}{2}\, \psi_{{1}}(t) ^{2}+\frac{1}{2}\,
\psi_{{2}}(t)  ^{2}+\psi_{{2}}(t)
x_{{1}}(t) \psi_{{3}}(t) +\frac{1}{2}\,
x_{1}(t)  ^{2} \psi_{{3}}(t)^{2}={\it const}
 \end{eqnarray*}
\normalsize
The last conservation law corresponds to the Hamiltonian.
This, and $\psi_3 = {\it const}$, are trivial conservation laws for the problem.
The missing first integral to solve the problem, $x_2 \psi_3 + \psi_1$,
was obtained in \cite{Taimanov97}.
\end{example}


\begin{example}[Cartan -- $\mathbf{(2,3,5)}$ problem]
\label{ex:Pr:235}
(30'34'')
The Cartan problem with growth vector $(2,3,5)$
can be posed in the following way:
\begin{eqnarray*}
&&\frac{1}{2}\,\int_a^b \left({u_1(t)}^2+{u_2(t)}^{2}\right)\mathrm{d}t \longrightarrow \min,\\
&&\left\{\begin{array}{l}
\dot x_1(t)=u_1(t)\, ,\\
\dot x_2(t)=u_2(t)\, ,\\
\dot x_3(t)=u_2(t) x_1(t)\, ,\\
\dot x_4(t)=\frac{1}{2}\,u_2(t) {x_1(t)}^2\, ,\\
\dot x_5(t)=u_2(t) x_1(t) x_2(t)\, .
\end{array}\right.
\end{eqnarray*}
The integrability of the problem was recently established in \cite{Sachkov04}.
This is possible with five independent conservation laws. They can easily be
determined with our \maple package.
\small
\begin{verbatim}
> L:=1/2*(u[1]^2+u[2]^2);
> phi:=[u[1], u[2], u[2]*x[1], 1/2*u[2]*x[1]^2, u[2]*x[1]*x[2]];
\end{verbatim}
\begin{eqnarray*}
L&:=&\frac{1}{2}\,{u_1}^2+\frac{1}{2}\,{u_2}^2\\
\varphi &:=&[u_1, u_2, u_2 x_1, \frac{1}{2}\,u_2 {x_1}^2, u_2 x_1 x_2]\\
\end{eqnarray*}
\normalsize
\small
\begin{verbatim}
> Symmetry(L, phi, t, [x[1],x[2],x[3],x[4],x[5]], [u[1],u[2]]);
\end{verbatim}
\begin{multline*}
\Biggl\{T=\frac{2}{3}\,C_{{1}}t+C_{{4}},\,
X_{{1}}=\frac{1}{3}\,C_{{1}}x_{{1}},\,
X_{{2}}=\frac{1}{3}\,C_{{1}}x_{{2}}+C_{{2}},\,
X_{{3}}=\frac{2}{3}\,C_{{1}}x_{{3}}+C_{{6}},\\
X_{{4}}=C_{{1}}x_{{4}}+C_{{5}},\,
X_{{5}}=C_{{2}}x_{{3}}+C_{{1}}x_{{5}}+C_{{3}},\,
U_{{1}}=-\frac{1}{3}\,C_{{1}}u_{{1}},\,
U_{{2}}=-\frac{1}{3}\,C_{{1}}u_{{2}},\\
\Psi_{{1}}=-\frac{1}{3}\,C_{{1}}\psi_{{1}},\,
\Psi_{{2}}=-\frac{1}{3}\,\psi_{{2}}C_{{1}},\,
\Psi_{{3}}=-\frac{2}{3}\,C_{{1}}\psi_{{3}}-C_{{2}}\psi_{{5}},\,
\Psi_{{4}}=-C_{{1}}\psi_{{4}},\,
\Psi_{{5}}=-C_{{1}}\psi_{{5}}\Biggr\}
\end{multline*}
\normalsize
\small
\begin{verbatim}
> CL:=Noether(L, phi, t, [x[1],x[2],x[3],x[4],x[5]], [u[1],u[2]], %);
\end{verbatim}
\begin{multline*}
 CL:=\frac{1}{3}\,C_{{1}}x_{{1}}(t) \psi_{{1}}(t) +
 \left( \frac{1}{3}\,C_{{1}}x_{{2}}(t) +C_{{2}} \right) \psi_{{2}}(t)
 + \left( \frac{2}{3}\,C_{{1}}x_{{3}}(t) +C_{{6}} \right) \psi_{{3}}(t)
+ \left( C_{{1}}x_{{4}}(t) + C_{{5}} \right) \psi_{{4}}(t) \\
 + \left( C_{{2}}x_{{3}}(t) +C_{{1}}x_{{5}} \left( t
 \right) +C_{{3}} \right) \psi_{{5}}(t)
 - \Biggl( \psi_{{0}} \left( \frac{1}{2}\,
 \left( u_{{1}}(t)  \right) ^{2}+\frac{1}{2}\,
 \left( u_{{2}}(t)  \right) ^{2} \right)
 +u_{{1}}(t) \psi_{{1}}(t) +u_{{2}} \left( t \right) \psi_{{2}}(t)\\
 +u_{{2}}(t) x_{{1}}(t) \psi_{{3}}(t)
 +\frac{1}{2}\,u_{{2}} \left( t
 \right)  \left( x_{{1}}(t)  \right) ^{2}\psi_{{4}}(t)
 +u_{{2}}(t) x_{{1}}(t) x_{{2}}(t) \psi_{{5}}(t)  \Biggr)
\left( \frac{2}{3}\,C_{{1}}t+C_{{4}} \right) ={\it const}
\end{multline*}
\normalsize
The Hamiltonian is given by
\small
\begin{verbatim}
> H:=-L+Vector[row]([psi[1](t), psi[2](t), psi[3](t), psi[4](t),
                                                         psi[5](t)]).Vector(phi);
\end{verbatim}
\[
H:=-\frac{1}{2}\,{u_{{1}}}^{2}-\frac{1}{2}\,{u_{{2}}}^{2}+u_{{1}}\psi_{{1}} \left( t
 \right) +u_{{2}}\psi_{{2}}(t) +u_{{2}}x_{{1}}\psi_{{3}}(t)
 +\frac{1}{2}\,u_{{2}}{x_{{1}}}^{2}\psi_{{4}} \left( t
 \right) +u_{{2}}x_{{1}}x_{{2}}\psi_{{5}}(t)
\]
\normalsize
and the extremal controls are obtained through the stationary condition \eqref{eq:stat:cond}.
\small
\begin{verbatim}
> solve({diff(H,u[1])=0, diff(H,u[2])=0}, {u[1],u[2]}):
> subs(x[1]=x[1](t),x[2]=x[2](t), u[1]=u[1](t),u[2]=u[2](t), %);
\end{verbatim}
\[
\left\{ u_{{1}} \left( t \right) =\psi_{{1}} \left( t \right) ,u_{{2}
} \left( t \right) =x_{{1}} \left( t \right) x_{{2}} \left( t \right)
\psi_{{5}} \left( t \right) +\psi_{{2}} \left( t \right) +x_{{1}}
 \left( t \right) \psi_{{3}} \left( t \right) +\frac{1}{2}\, x_{{1}}
 \left( t \right)  ^{2}\psi_{{4}} \left( t \right)  \right\}
\]
\normalsize

\small
\begin{verbatim}
> CL:=subs(psi[0]=-1, %, CL);
\end{verbatim}
\begin{multline*}
CL:=\frac{1}{3}\,C_{{1}}x_{{1}}(t) \psi_{{1}}(t) +
 \left( \frac{1}{3}\,C_{{1}}x_{{2}}(t) +C_{{2}} \right) \psi_{{2}}(t)
 + \left( \frac{2}{3}\,C_{{1}}x_{{3}}(t) +C_{{6}} \right) \psi_{{3}}(t)\\
+ \left( C_{{1}}x_{{4}}(t) + C_{{5}} \right) \psi_{{4}}(t)
 +  \left( C_{{2}}x_{{3}}(t) +C_{{1}}x_{{5}} \left( t
 \right) +C_{{3}} \right) \psi_{{5}}(t)
 - \Biggr( \frac{1}{2}\,\psi_{{1}}(t) ^{2}\\
 -\frac{1}{2}\, \left( x_{{1}}(t) x_{{2}}(t) \psi_{{5}} \left( t
 \right) +\psi_{{2}}(t) +x_{{1}}(t) \psi_{{3}}(t) +\frac{1}{2}\,  x_{{1}}(t)
^{2}\psi_{{4}}(t)  \right) ^{2}\\
+\left( x_{{1}} \left( t
 \right) x_{{2}}(t) \psi_{{5}}(t) +\psi_{{2}}(t)
 + x_{{1}}(t) \psi_{{3}} \left( t
 \right) +\frac{1}{2}\, x_{{1}}(t)  ^{2}\psi_{{4}}(t)  \right) \psi_{{2}}(t)\\
+\left( x_{{1}}(t) x_{{2}}(t) \psi_{{5}} \left( t
 \right) +\psi_{{2}}(t) +x_{{1}}(t) \psi_{{3}}(t)
 + \frac{1}{2}\, x_{{1}}(t)^{2}\psi_{{4}}(t)  \right) x_{{1}}(t)\psi_{{3}}(t)\\
+\frac{1}{2}\, \left( x_{{1}}(t) x_{{2}}(t) \psi_{{5}}(t) +\psi_{{2}}\left(t\right)
+ x_{{1}}(t) \psi_{{3}}(t) +\frac{1}{2}\,
 \left( x_{{1}}(t)  \right) ^{2}\psi_{{4}} \left( t
 \right)  \right)  x_{{1}}(t)^{2}\psi_{{4}}(t)\\
+\left( x_{{1}}(t) x_{{2}}(t) \psi_{{5}}(t) +\psi_{{2}} \left( t
 \right) +x_{{1}}(t) \psi_{{3}}(t) +\frac{1}{2}\,
 x_{{1}}(t) ^{2}\psi_{{4}} \left( t
 \right)  \right) x_{{1}}(t) x_{{2}}(t)
\psi_{{5}}(t)  \Biggl)
\left( \frac{2}{3}\,C_{{1}}t+C_{{4}}\right) ={\it const}
\end{multline*}
\normalsize
The five conservation laws we are looking for, are easily obtained
(the last one corresponds to the Hamiltonian):
\small
\begin{verbatim}
> subs(C[6]=1,seq(C[i]=0,i=1..6), CL);
> subs(C[5]=1,seq(C[i]=0,i=1..6), CL);
> subs(C[3]=1,seq(C[i]=0,i=1..6), CL);
> subs(C[2]=1,seq(C[i]=0,i=1..6), CL);
> simplify(subs(C[4]=-1,seq(C[i]=0,i=1..6), CL));
\end{verbatim}
\begin{eqnarray*}
&&\psi_{{3}}(t) ={\it const}\\
&&\psi_{{4}}(t) ={\it const}\\
&&\psi_{{5}}(t) ={\it const}\\
&&\psi_{{2}}(t) +\psi_{{5}}(t) x_{{3}}(t) ={\it const}\\
&&\frac{1}{2}\,  x_{{1}}(t) ^{2} x_{{2}}(t)  ^{2} \psi_{{5}}(t)^{2}
+x_{{1}}(t) x_{{2}}(t) \psi_{{5}}(t) \psi_{{2}}(t)
+ x_{{1}}(t)^{2} x_{{2}}(t)\psi_{{5}}(t) \psi_{{3}}(t)
+\frac{1}{2}\, x_{{1}}(t)^{3}x_{{2}}(t) \psi_{{5}}(t) \psi_{{4}}(t)\\
&&  +\frac{1}{2}\, x_{{1}}(t)^{2} \psi_{{3}}(t)^{2}
+\frac{1}{8}\,  x_{{1}}(t)^{4}\psi_{{4}}(t)^{2}
+\frac{1}{2}\,  \psi_{{1}}(t)^{2}+\psi_{{2}}(t) x_{{1}}(t) \psi_{{3}}(t)
+\frac{1}{2}\,\psi_{{2}}\left(t\right) x_{{1}}(t)^{2}\psi_{{4}}(t)\\
&& +\frac{1}{2}\, x_{{1}}(t)^{3}
\psi_{{3}}(t) \psi_{{4}}(t) +\frac{1}{2}\,\psi_{{2}}(t)^{2}={\it const}
\end{eqnarray*}
\normalsize
One can say that for the Cartan $(2,3,5)$ problem we have four
trivial first integrals: the Hamiltonian $H$; and the multipliers
$\psi_3$, $\psi_4$, $\psi_5$. Together with the non-trivial
integral $x_3 \psi_5 + \psi_2$, the problem becomes completely
integrable (see \cite{Sachkov04}).
\end{example}


\section{Conservation Laws in the Calculus of Variations}
\label{sec:CV}

Let us consider the classical problem of the calculus
of variations with higher-order derivatives:
to minimize an integral functional
\begin{equation}
\label{eq:funcional}
J[\mathbf{x}(\cdot)] = \int_{a}^{b}
L(t,\mathbf{x}(t),\dot{\mathbf{x}}(t),\ldots,\mathbf{x}^{(r)}(t))
\,\textrm{d}t \, ,
\end{equation}
subject to certain boundary conditions, and
where the Lagrangian $L$ depends on the independent
variable $t \in \mathbb{R}$, on $n$ dependent variables
$\mathbf{x}(t)=[ x_1(t),\, x_2(t),\,\cdots,\,x_n(t)]^\textsuperscript{T} \in \mathbb{R}^n$,
and its $r$ first derivatives ($\dot{\mathbf{x}}(t)\equiv \mathbf{x}^{(1)}(t)$):
\begin{equation*}
\mathbf{x}^{(i)}(t)=\left[ \frac{\mathrm{d}^i x_1(t)}{\mathrm{d} t^i}, \,
\frac{\mathrm{d}^ix_2(t)}{\mathrm{d}t^i},\, \cdots,\,
\frac{\mathrm{d}^ix_n(t)}{\mathrm{d}t^i}  \right]^\textsuperscript{T}
\in \mathbb{R}^n \, , \quad i=1,\,\ldots,\, r\, .
\end{equation*}
It is well known that the problems of the calculus of variations
are a particular case of the optimal control problem
(\ref{eq:funcionalCO})-(\ref{eq:sistCont}). The standard
technique to write the problem of minimizing \eqref{eq:funcional}
as an optimal control problem consists to introduce new functions
\[
\mathbf{x}^0=\mathbf{x},\; \mathbf{x}^1=\dot{\mathbf{x}},\;
\mathbf{x}^2=\mathbf{x}^{(2)},\; \ldots,\;
\mathbf{x}^{r-1}=\mathbf{x}^{(r-1)},\; \textrm{ and } \mathbf{u}=\mathbf{x}^{(r)} \, .
\]
With this notation at hand, the equivalent optimal control problem
has $rn$ state variables ($x_i^j$, $i=1,\ldots,n$, $j=0,\ldots,r-1$),
and $n$ controls ($\mathbf{u}=\mathbf{x}^{(r)}$):
\begin{gather*}
J[\mathbf{x}(\cdot),\mathbf{u}(\cdot)] = \int_{a}^{b}
L(t,\mathbf{x}(t),\mathbf{u}(t)) \,\textrm{d}t \longrightarrow \min \, , \\
\begin{cases}
 \dot{\mathbf{x}}^0(t)&=\mathbf{x}^1(t)\, ,\\
 \dot{\mathbf{x}}^1(t)&=\mathbf{x}^2(t)\, ,\\
 &\vdots\\
 \dot{\mathbf{x}}^{r-2}(t)&=\mathbf{x}^{r-1}(t)\, ,\\
 \dot{\mathbf{x}}^{r-1}(t)&=\mathbf{u}(t)\, .
\end{cases}
\end{gather*}
Since any problem of the calculus of variations can
always be rewritten as an optimal control problem,
we can also apply our \maple package to obtain variational
symmetries and conservation laws in the classical context
of the calculus of variations, and thus recovering
the previous investigations of the authors \cite{gouv04}.
We recall that for the problems of the calculus of variations
there are no abnormal extremals (one can always choose $\psi_0 = -1$).
Follow some examples.


\begin{example}
\label{ex:Ex11}
(0'08'')
We begin with a very simple problem of the calculus of variations,
where the Lagrangian depends only on one dependent variable
($n = 1$), and where there are no derivatives of higher order than the first one ($r = 1$):
$L(t,x,\dot{x})=t\dot{x}^2$. According with the above mentioned technique
of rewriting the problem as an optimal control problem,
we write the following definitions in \mapleSE :
\small
\begin{verbatim}
> L:=t*v^2; u:=v; phi:=u;
\end{verbatim}
\begin{eqnarray*}
L &:=& t{v}^{2}\\
u &:=& v\\
\varphi &:=& v
\end{eqnarray*}
\normalsize
Our procedure \emph{Symmetry} determine the general infinitesimal generators
which define the family of symmetries for the problem of the calculus of variations
under consideration:
\small
\begin{verbatim}
> Symmetry(L,phi,t,x,u);
\end{verbatim}
\[
 \left\{ T=t C_{1},\, X=C_{2},\, U=-v C_{1},\, \Psi=0 \right\}
\]
\normalsize
The conservation laws corresponding to the computed
symmetries, are obtained with the procedure \emph{Noether}:
\small
\begin{verbatim}
> Noether(L,phi,t,x,u,%);
\end{verbatim}
\[
C_{2}\psi(t) - \left( \psi_{0}t \left( v(t)
\right) ^{2}+\psi(t) v(t)  \right) tC_{1}
\mbox{}={\it const}
\]
\normalsize
Going back to the original notation ($v = \dot{x}$),
we can write:
\small
\begin{verbatim}
> CL:=subs(psi[0]=-1,v(t)=diff(x(t),t),%);
\end{verbatim}
\[
{\it CL}:=C_{2}\psi(t) + \left(t
\left( {\frac {d}{dt}}x(t)  \right) ^{2}-\psi(t)
{\frac {d}{dt}}x(t)  \right) tC_{1}
\mbox{}={\it const}.
\]
\normalsize
In this case one can easily use the definition
of first integral (a function that is preserved along
the extremals of the problem), to verify the validity
of the obtained expression. For that we compute the pair $(x(t),\psi(t))$ that satisfies
the adjoint system (\ref{eq:sistAdj}) and the maximality condition
(\ref{eq:condMax}) of the Pontryagin maximum principle
(Theorem~\ref{thm:pMaxPont}).
\small
\begin{verbatim}
> H:=-L+psi(t)*phi;
\end{verbatim}
\[
H:=-t{v}^{2}+\psi(t) v
\]
\normalsize
\small
\begin{verbatim}
> {diff(H,u)=0, diff(psi(t),t)=-diff(H,x)};
\end{verbatim}
\[
 \left\{ {\frac {d}{dt}}\psi(t) =0, -2\,tv+\psi
(t) =0 \right\}
\]
\normalsize
After substituting $v = \dot{x}$, we
obtain the extremals by solving the above system
of differential equations:
\small
\begin{verbatim}
> subs(v=diff(x(t),t),%);
\end{verbatim}
\[
 \left\{ {\frac {d}{dt}}\psi(t) =0, -2\,t{\frac {d}{dt}}x
(t) +\psi(t) =0 \right\}
\]
\normalsize
\small
\begin{verbatim}
> dsolve(%);
\end{verbatim}
\[
 \left\{ \psi(t) ={\it K_2},x(t)
 =\frac{1}{2}\,{\it K_2}\,\ln(t)+{\it K_1}
 \right\}
\]
\normalsize
Expression for $x(t)$ coincides with the Euler-Lagrange extremal
(\cite[Example~5.1]{gouv04}). Substituting the extremals in the
conservation law one obtains, as expected, a true proposition:
\small
\begin{verbatim}
> expand(subs(%,CL));
\end{verbatim}
\[
C_{2}{\it K_2}-\frac{1}{4}\,C_{1}{{\it K_2}}^{2}={\it const}
\]
\normalsize
Substituting only $\psi(t)$, one can get the family of conservation laws
in the notation of the calculus of variations:
\small
\begin{verbatim}
> expand(subs(psi(t)=K[2],CL));
\end{verbatim}
\[
C_{2}{\it K_2}+{t}^{2}C_{1} \left( {\frac {d}{dt}}x
(t)  \right) ^{2}-tC_{1}{\it K_2}\,{\frac {d}{dt}}x
(t)={\it const}
\]
\normalsize
\small
\begin{verbatim}
> subs(C[2]*K[2]=0,C[1]=-1,%);
\end{verbatim}
\[
-{t}^{2} \left( {\frac {d}{dt}}x(t)  \right) ^{
2}+t K_2{\frac {d}{dt}}x(t) ={\it const}
\]
\normalsize
\end{example}


\begin{example}[Kepler's problem]
\label{ex:Ex12}
(0'17'')
We now obtain the conservation laws
for Kepler's problem -- see \cite[p.~217]{BRU}.
In this case the Lagrangian depends on
two dependent variables ($n = 2$), without
derivatives of higher-order ($r = 1$):
\[
L(t,\mathbf{q},\dot{\mathbf{q}})=\frac{m}{2}\left(\dot{q}_1^2
+\dot{q}_2^2\right)+\frac{K}{\sqrt{q_1^2+q_2^2}} \, .
\]
The family of conservation laws for the problem
is easily obtained with our \maple package:
\small
\begin{verbatim}
> L:= m/2*(v[1]^2+v[2]^2)+K/sqrt(q[1]^2+q[2]^2); x:=[q[1],q[2]]; u:=[v[1],v[2]];
  phi:=[v[1],v[2]];
\end{verbatim}
\begin{eqnarray*}
L &:=& \frac{1}{2}\,m \left( {v_{1}}^{2}+{v_{2}}^{2} \right) +{\frac {K}{\sqrt
{{q_{1}}^{2}+{q_{2}}^{2}}}}\\
x &:=& [q_{1},q_{2}]\\
u &:=& [v_{1},v_{2}]\\
\varphi &:=& [v_{1},v_{2}]
\end{eqnarray*}
\begin{verbatim}
> Symmetry(L, phi, t, x, u);
\end{verbatim}
\begin{multline*}
\Biggl\{ T=C_{3},\,
X_{1}=-C_{1}q_{2},\, X_{2}=C_{1}q_{1},\,
U_{1}=-C_{2}v_{2}+{\frac { \left( C_{1}-C_{2}\right) \psi_{2}}{\psi_{0}m}},\,
U_{2}=C_{2}v_{1}-{\frac { \left( C_{1}-C_{2} \right) \psi_{1}}{\psi_{0}m}},\\
\Psi_{1}=-C_{1}\psi_{2},\, \Psi_{2}=C_{1}\psi_{1}\Biggr\}
\end{multline*}
\begin{verbatim}
> Noether(L, phi, t, x, u, %):
> CL:=subs(psi[0]=-1,%);
\end{verbatim}
\begin{multline*}
{\it CL}:=-C_{1}q_{2}(t) \psi_{1}(t)
+C_{1}q_{1}(t) \psi_{2}(t)\\
- \Biggl(-\frac{1}{2}\,m \left(  v_{1}(t)
 ^{2}+  v_{2}(t)  ^{2} \right)
-{\frac {K}{\sqrt { q_{1}(t) ^{2}+
q_{2}(t) ^{2}}}}
+\psi_{1}(t) v_{1}(t) +\psi_{2}
(t) v_{2}(t)  \Biggr) C_{3}={\it const}
\end{multline*}
\normalsize
To obtain the conservation laws in the format of the calculus
of variations, one needs to compute the Pontryagin multipliers
$(\psi_{1}(t),\psi_{2}(t))$,
\small
\begin{verbatim}
> H:=-L+Vector[row]([psi[1](t), psi[2](t)]).Vector(phi);
\end{verbatim}
\[
H:=-\frac{1}{2}\,m \left( {v_{1}}^{2}+{v_{2}}^{2} \right)
-{\frac {K}{\sqrt {{q_{1}}^{2}+{q_{2}}^{2}}}} +v_{1}\psi_{1}
(t)+v_{2}\psi_{2}(t)
\]
\normalsize
\small
\begin{verbatim}
> solve({diff(H,v[1])=0,diff(H,v[2])=0},{psi[1](t), psi[2](t)});
\end{verbatim}
\[
 \left\{ \psi_{1}(t) =mv_{1},\psi_{2}
(t) =mv_{2} \right\}
\]
\normalsize
and substitute the expressions, together with
$v_{1}(t) = \dot q_{1}(t)$ and $v_{2}(t) = \dot q_{2}(t)$:
\small
\begin{verbatim}
> expand(subs(%,v[1](t)=v[1],v[2](t)=v[2],v[1]=diff(q[1](t),t),
                v[2]=diff(q[2](t),t),CL));
\end{verbatim}
\begin{eqnarray*}
-C_{1}q_{2}(t) m{\frac {d}{dt}}q_{1}
(t)+C_{1}q_{1}(t) m
{\frac {d}{dt}}q_{2}(t)-\frac{1}{2}\,C_{3}m
\left( {\frac{d}{dt}}q_{1}(t)  \right) ^{2}-\frac{1}{2}\,C_{3}
m \left( {\frac {d}{dt}}q_{2}(t)  \right) ^{2}\\ \hbox{}
+{\frac {C_{3}K}{\sqrt { q_{1}(t)
 ^{2}+  q_{2}(t) ^{2}}}}
={\it const}
\end{eqnarray*}
\normalsize
This is the conservation law in \cite[Example~5.2]{gouv04}.
\end{example}


\begin{example}
\label{ex:Ex13}
(6'42'')
Let us see an example of the calculus of variations
whose Lagrangian depends on two functions ($n=2$) and higher-order
derivatives ($r=2$):
\[
L(t,\mathbf{x},\dot{\mathbf{x}},\ddot{\mathbf{x}})=\dot{x}_1^2+\ddot{x}_2^2 \, .
\]
We write the problem in the optimal control terminology, and make use
of our \maple procedure \emph{Symmetry} to compute the symmetries:
\small
\begin{verbatim}
> L:=v[1]^2+a[2]^2; xx:=[x[1],x[2],v[1],v[2]]; u:=[a[1],a[2]];
  phi:=[v[1],v[2],a[1],a[2]];
\end{verbatim}
\begin{eqnarray*}
L &:=& {v_{1}}^{2}+{a_{2}}^{2}\\
{\it xx} &:=& [x_{1},x_{2},v_{1},v_{2}]\\
u &:=& [a_{1},a_{2}]\\
\varphi &:=& [v_{1},v_{2},a_{1},a_{2}]
\end{eqnarray*}
\begin{verbatim}
> Symmetry(L, phi, t, xx, u);
\end{verbatim}
\begin{multline*}
\biggl\{T=2\,C_{3}t+C_{4},\,
X_{1}=C_{3}x_{1}+C_{5},\,
X_{2}=tC_{1}+3\,C_{3}x_{2}+C_{6},\,
X_{3}=-C_{3}v_{1},\,
X_{4}=C_{1}+C_{3}v_{2},\\
U_{1}=-3\,C_{3}a_{1}-2\,\psi_{0}C_{2}a_{2}-C_{2}\psi_{4},\,
U_{2}=-C_{3}a_{2}+C_{2}\psi_{3},\\
\Psi_{1}=-C_{3}\psi_{1},\,
\Psi_{2}=-3\,C_{3}\psi_{2},\,
\Psi_{3}=C_{3}\psi_{3},\,
\Psi_{4}=-C_{3}\psi_{4}\biggr\}
\end{multline*}
\normalsize
We choose, in the conservation law returned by our \maple procedure
\emph{Noether}, $\psi_0=-1$, and then go back to the
calculus of variations notation:
$v_1(t) = \dot x_1(t)$, $v_2(t) = \dot x_2(t)$,
$a_1(t) = \ddot x_1(t)$, and $a_2(t) = \ddot x_2(t)$:
\small
\begin{verbatim}
> Noether(L, phi, t, xx, u, %):
> CL:=subs(psi[0]=-1, v[1](t)=diff(x[1](t),t), v[2](t)=diff(x[2](t),t),
  a[1](t)=diff(x[1](t),t$2), a[2](t)=diff(x[2](t),t$2),%);\end{verbatim}
\begin{eqnarray*}
{\it CL}:= \left( C_{3}x_{1}(t) +C_{5} \right)
\psi_{1}(t) + \left( tC_{1}+3\,C_{3}x_{2}(t)
+C_{6} \right) \psi_{2}(t)
\mbox{}-C_{3} \left( {\frac {d}{dt}}x_{1}(t)  \right)
\psi_{3}(t) \\
+\left( C_{1}+C_{3}{\frac {d}{dt}}x_{2}(t)  \right)
\psi_{4}(t) - \biggl( - \left( {\frac{d}{dt}}
x_{1}(t)  \right) ^{2}- \left( {\frac {d^2}{d t^2}}x_{2}
(t)  \right) ^{2} +\psi_{1}(t)
{\frac {d}{dt}}x_{1}(t)\\
\mbox{}+\psi_{2}(t) {\frac {d}{dt}}x_{2}(t)
+\psi_{3}(t) {\frac {d^2}{d t^2}}x_{1}(t)
+\psi_{4}(t) {\frac {d^2}{d t^2 }}x_{2}(t)
\biggr)  \left( 2\,C_{3}t+C_{4} \right)
={\it const}
\end{eqnarray*}
\normalsize
Similarly to Example~\ref{ex:Ex11}, we can also compute through \maple
the extremals and verify, by definition, the validity
of the obtained family of conservation laws.
\small
\begin{verbatim}
> vpsi:=Vector[row]([psi[1](t), psi[2](t), psi[3](t), psi[4](t)]):
> H:=-L+vpsi.Vector(phi);
\end{verbatim}
\[
H:=-{v_{1}}^{2}-{a_{2}}^{2} +v_{1}
\psi_{1}(t) +v_{2}\psi_{2}(t) +a_{1}
\psi_{3}(t)+a_{2}\psi_{4}(t)
\]
\normalsize
\small
\begin{verbatim}
> {diff(H,u[1])=0, diff(H,u[2])=0, diff(vpsi[1],t)=-diff(H,xx[1]),
   diff(vpsi[2],t)=-diff(H,xx[2]), diff(vpsi[3],t)=-diff(H,xx[3]),
   diff(vpsi[4],t)=-diff(H,xx[4])}:
> subs(v[1]=diff(x[1](t),t), a[2]=diff(x[2](t),t$2), %);
\end{verbatim}
\begin{eqnarray*}
\left\{ {\frac {d}{dt}}\psi_{3}(t) =2\,
{\frac {d}{dt}}x_{1}(t) -\psi_{1}(t) ,
{\frac {d}{dt}}\psi_{2}(t) =0,{\frac {d}{dt}}\psi_{4}
(t) =-\psi_{2}(t) ,\psi_{3}(t)=0,
{\frac {d}{dt}}\psi_{1}(t) =0,\right.\\
\left.
-2\,{\frac {d^2}{d t^2 }}x_{2}(t)
\mbox{}+\psi_{4}(t) =0 \right\}
\end{eqnarray*}
\normalsize
Solving the above system of equations, that result from the maximality
condition and adjoint system,
\small
\begin{verbatim}
> dsolve(%);
\end{verbatim}
\begin{eqnarray*}
\Bigl\{ x_{1}(t) =\frac{1}{2}\,{\it K_6}\,t +{\it K_4}, \,
x_{2}(t) =-\frac{1}{12}\,{\it K_5}\,{t}^{3}+\frac{1}{4}\,{\it K_3}\,{t}^{2}
+{\it K_1}\,t+{\it K_2},\psi_{3}(t) =0,\\
\psi_{1}(t) ={\it K_6},\,
\psi_{4}(t)=-{\it K_5}\,t+{\it K_3},\, \psi_{2}(t) ={\it K_5}\Bigr\}
\end{eqnarray*}
\normalsize
we obtain the extremals (the extremal state trajectories $x_1(t)$ and $x_2(t)$
are the same as the ones obtained in \cite[Example~5.3]{gouv04},
by solving the Euler-Lagrange necessary optimality condition)
that, substituted in the conservation law,
\small
\begin{verbatim}
> expand(subs(%,CL));
\end{verbatim}
\begin{eqnarray*}
{\it K_6}\,C_{3}{\it K_4}+{\it K_6}\,C_{5}+3\,{\it K_5}\,C_{3}
{\it K_2}+{\it K_5}\,C_{6}+C_{1}{\it K_3}+C_{3}{\it K_1}\,
{\it K_3}-\frac{1}{4}\,{{\it K_6}}^{2}C_{4}\\
-\frac{1}{4}\,{{\it K_3}}^{2}C_{4}-{\it K_5}\,{\it K_1}\,C_{4}={\it const}
\end{eqnarray*}
\normalsize
conduces to a true proposition (constant equal constant).
Finally, substituting only the Pontryagin multipliers,
\small
\begin{verbatim}
> subs({psi[1](t)=K[6], psi[3](t)=0, psi[4](t)=-K[5]*t+K[3], psi[2](t)=K[5]}, CL);
\end{verbatim}
\begin{eqnarray*}
\left( C_{3}x_{1}(t) +C_{5} \right) {\it K_6}
+ \left( C_{1}t+3\,C_{3}x_{2}(t)+C_{6} \right){\it K_5}
+ \left( C_{1}+C_{3}{\frac {d}{dt}}x_{2}(t)  \right)
\left( -{\it K_5}\,t+{\it K_3} \right)\\ \hbox{}
- \biggl( - \left( {\frac {d}{dt}}x_{1}(t)
\right) ^{2}-\left( {\frac {d^2}{d t^2}}x_{2}(t)  \right) ^{2}
+{\it K_6}\,{\frac {d}{dt}}x_{1}(t) +{\it K_5}\,
{\frac {d}{dt}}x_{2}(t)\\ \hbox{}
+ \left( -{\it K_5}\,t+{\it K_3} \right) {\frac {d^2}{d t^2}}x_{2}
(t)  \biggr)  \left( 2\,C_{3}t+C_{4} \right) ={\it const}
\end{eqnarray*}
\normalsize
the conservation law takes the form of a differential equation of less order
than the one obtained in \cite[Example~5.3]{gouv04} (the Hamiltonian approach
is here more suitable than the Lagrangian one).
\end{example}


\begin{example}[Emden-Fowler]
\label{ex:Ex14}
(0'01'')
Given the variational problem of Emden-Fowler \cite[p. 220]{BRU},
defined by the Lagrangian
\small
\begin{verbatim}
> L:= t^2/2*(v^2-x^6/3);
\end{verbatim}
\[
L:=\frac{1}{2}\,{t}^{2} \left( {v}^{2}-\frac{1}{3}\,{x}^{6} \right)
\]
\normalsize
we are interested to find, following our methodology, the conservation laws for the problem.
\small
\begin{verbatim}
> Symmetry(L, v, t, x, v);
\end{verbatim}
\[
 \left\{ T=-2\,tC_{1},\, X=C_{1}x,\, U=3\,C_{1}v,\, \Psi=-\psi\,C_{1} \right\}
\]
\begin{verbatim}
> Noether(L,v,t,x,v,%):
> CL:=subs(psi[0]=-1,%);
\end{verbatim}
\[
{\it CL}:=C_{1}x(t) \psi(t) +2
\left( -\frac{1}{2}\,{t}^{2} \left( v(t)
^{2}-\frac{1}{3}\,  x(t)^{6} \right)
+\psi(t) v(t)  \right) tC_{1}={\it const}
\]
\normalsize
The expression for the $\psi(t)$ comes from
the stationary condition.
 \small
\begin{verbatim}
> H:=-L+psi(t)*v;
\end{verbatim}
\[
H:=-\frac{1}{2}\,{t}^{2} \left( {v}^{2}-\frac{1}{3}\,{x}^{6} \right)
+\psi(t) v
\]
\normalsize
\small
\begin{verbatim}
> solve(diff(H,v)=0,{psi(t)});
\end{verbatim}
\[
 \left\{ \psi(t) ={t}^{2}v \right\}
\]
\normalsize
\small
\begin{verbatim}
> subs(%,v(t)=diff(x(t),t),v=diff(x(t),t),CL): expand(%);
\end{verbatim}
\[
C_{1}x(t) {t}^{2}{\frac {d}{dt}}x(t)
+{t}^{3}C_{1} \left( {\frac {d}{dt}}x(t)
\right) ^{2}+\frac{1}{3}\,{t}^{3}C_{1} x(t)
^{6}={\it const}
\]
\normalsize
Fixing $C_1 = 3$,
\small
\begin{verbatim}
> subs(C[1]=3,%);
\end{verbatim}
\[
3\,x(t) {t}^{2}{\frac {d}{dt}}x(t) +3\,{t}^{3}
\left( {\frac {d}{dt}}x(t)  \right) ^{2}+{t}^{3}
x(t)  ^{6}={\it const}
\]
\normalsize
we obtain the same conservation law as the one obtained in \cite[Example~5.4]{gouv04},
with the methods of the calculus of variations.
\end{example}


\begin{example}[Thomas-Fermi]
\label{ex:Ex15}
(0'01'')
We consider the problem of Thomas-Fermi \cite[p.~220]{BRU},
showing an example of a problem of the calculus of variations which does
not admit variational symmetries.
\small
\begin{verbatim}
> L:=1/2*v^2+2/5*x^(5/2)/sqrt(t);
\end{verbatim}
\[
L:=\frac{1}{2}\,{v}^{2}+\frac{2}{5}\,{\frac {{x}^{\frac{5}{2}}}{\sqrt {t}}}
\]
\begin{verbatim}
> Symmetry(L, v, t, x, v);
\end{verbatim}
\[
\left\{ T=0,\, X=0,\, U=0,\, \Psi=0 \right\}
\]
\normalsize
Our \maple function \emph{Symmetry} returns, in this case, vanishing generators.
As explained in \S\ref{sec:CCL}, this means that the problem
does not admit symmetries.
\end{example}


\begin{example}[Damped Harmonic Oscillator]
\label{ex:Ex16}
(0'02'')
We consider a harmonic oscillator
with restoring force $-kx$,
emersed in a liquid in such a way that the motion of the mass $m$
is damped by a force proportional to its velocity.
Using Newton's second law one obtains, as the equation of motion,
the Euler-Lagrange differential equation associated with
the following Lagrangian \cite[pp. 432--434]{LOG}:
\small
\begin{verbatim}
> L:=1/2*(m*v^2-k*x^2)*exp((a/m)*t);
\end{verbatim}
\[
L:=\frac{1}{2}\, \left( m{v}^{2}-k{x}^{2} \right) {e^{{\frac {at}{m}}}}
\]
\normalsize
In order to find the conservation laws, we first obtain,
as usual, the generators which define the symmetries of the problem.
\small
\begin{verbatim}
> simplify(Symmetry(L,v,t,x,v));
\end{verbatim}
\[
\left\{T=-{\frac {2\,mC_{1}}{a}},\,
X=C_{1}x,\, U=C_{1}v,\, \Psi=-\psi\,C_{1} \right\}
\]
\normalsize
\small
\begin{verbatim}
> Noether(L,v,t,x,v,%):
> CL:=subs(psi[0]=-1,%);
\end{verbatim}
\[
{\it CL}:=C_{1}x(t) \psi(t) + 2 \left( -\frac{1}{2}\,
 \left( m \,v(t)  ^{2}-k \, x(t)
^{2} \right) {e^{{\frac {at}{m}}}}+\psi
(t) v(t)  \right) m C_{1}{a}^{-1} ={\it const}
\]
\normalsize
The value for $\psi(t)$ is easily determined,
and we can write the obtained family of conservation
laws in the language of the calculus of variations.
\small
\begin{verbatim}
> H:=-L+psi(t)*v;
\end{verbatim}
\[
H:=-\frac{1}{2}\,\left( m{v}^{2}-k{x}^{2} \right)
{e^{{\frac {at}{m}}}}+\psi(t) v
\]
\normalsize
\small
\begin{verbatim}
> solve(diff(H,v)=0,{psi(t)});
\end{verbatim}
\[
 \left\{ \psi(t) =m\,v\,{e^{{\frac {at}{m}}}} \right\}
\]
\normalsize
\small
\begin{verbatim}
> simplify(subs(%,v(t)=diff(x(t),t),v=diff(x(t),t),CL));
\end{verbatim}
\[
C_{1}m{e^{{\frac {at}{m}}}} \left( x(t)  \left(
{\frac {d}{dt}}x(t)  \right) a+m \left( {\frac {d}{dt}}x
(t)  \right) ^{2}+k x(t) ^{2}
\right){a}^{-1}={\it const}
\]
\normalsize
Choosing an appropriate value to the constant $C_1$
\small
\begin{verbatim}
> subs(C[1]=-a/(2*m),%);
\end{verbatim}
\[
-\frac{1}{2}\,{e^{{\frac {at}{m}}}} \left( x(t)  \left( {\frac {d}{dt}}x
(t)  \right) a+m \left( {\frac {d}{dt}}x(t)\right)
^{2}+k x(t) ^{2} \right) ={\it const}
\]
\normalsize
we obtain the conservation law in \cite[Ch.~7, Example~1.10]{LOG}.
\end{example}


\section{The \maple package}
\label{sec:ProcMaple}

The procedures \emph{Symmetry} and \emph{Noether}, described in
the previous sections, have been implemented for the computer
algebra system \maple (version~9.5).
\begin{description}
\item[Symmetry] computes the infinitesimal generators which define
the symmetries of the optimal control problem specified in the
input. As explained in sections \ref{sec:CL} and \ref{sec:CCL},
this procedure involves the solution of a system of partial
differential equations. We have used the \maple solver
\emph{pdsolve}, trying to separate the variables by sum. \item
Output:
\begin{itemize}
\item[-] set of infinitesimal generators.
\end{itemize}
\item Syntax:
\begin{itemize}
\item[-] Symmetry(L, $\varphi$, t, x, u, [\texttt{all}])
\end{itemize}
\item Input:
\begin{itemize}
\item[L -] expression of the Lagrangian;
\item[$\varphi$ -] expression or list of expressions of the velocity vector
$\boldsymbol{\varphi}$ which defines the control system;
\item[t -] name of the independent variable;
\item[x -] name or list of names of the state variables;
\item[u -] name or list of names of the control variables;
\item[\texttt{all} -] This is an optional parameter. When \emph{all} is given
in the last argument of the procedure \emph{Symmetry},
the output presents all the constants given by the \maple command
\emph{pdsolve}. By default, that is, without optional \emph{all},
we eliminate redundant constants. This is done by our \maple
procedure \emph{reduzConst}. This is a technical routine,
and thus not provided here. Essentially, the procedure
transforms in one constant each sum of constants not repeated
in the conservation law. The interested reader can find the \maple file
with its definition, together with the \emph{Symmetry} and \emph{Noether}
code, that constitute our \maple package, at \textsf{http://www.mat.ua.pt/delfim/maple.htm}.
\end{itemize}
\end{description}
\small
\begin{verbatim}
Symmetry:=proc(L::algebraic, phi::{algebraic, list(algebraic)}, t::name,
                                    x0::{name,list(name)}, u0::{name,list(name)})
  local n,m, xx, i, vX, vPSI, vU, vv, lpsi, H, eqd, syseqd, sol, conjGerad, lphi;
  unprotect(Psi); unassign('T'); unassign('X'); unassign('U'); unassign('Psi');
  unassign('psi');
  n:=nops(x0); m:=nops(u0);
  if n>1 then lphi:=phi;lpsi:=[seq(psi[i],i=1..n)];
  else lphi:=[phi]; lpsi:=[psi]; fi;
  xx:=op(x0),op(u0),op(lpsi); vv:=Vector([seq(v||i,i=1..2*n+m)]);
  if n>1 then vX:=Vector([seq(X[i](t,xx), i=1..n)]);
  else vX:=Vector([X(t,xx)]); fi;
  if n>1 then vPSI:=Vector([seq(PSI[i](t,xx),i=1..n)]);
  else vPSI:=Vector([PSI(t,xx)]); fi;
  if m>1 then vU:=Vector([seq(U[i](t,xx), i=1..m)]);
  else vU:=Vector([U(t,xx)]); fi;
  H:=psi[0]*L+Vector[row](lphi).Vector(lpsi);
  eqd:=diff(H,t)*T(t,xx) +Vector[row]([seq(diff(H,i),i=x0)]).vX+Vector[row]([seq(
       diff(H,i),i=u0)]).vU+Vector[row]([seq(diff(H,xx[i]),i=n+m+1..n+m+n)]).vPSI
       -LinearAlgebra[Transpose](vPSI).vv[1..n]-Vector[row](lpsi).(map(diff,vX,t)
                       +Matrix([seq(map(diff,vX,i),i=xx)]).vv)+H*(diff(T(t,xx),t)
                                   +Vector[row]([seq(diff(T(t,xx),i),i=xx)]).vv);
  eqd:=expand(eqd);  eqd:=collect(eqd, convert(vv,'list'), distributed);
  syseqd:={coeffs(eqd, convert(vv,'list'))}:
  conjGerad:={T(t,xx)}union convert(vX,'set') union convert(vU,'set')
                                                       union convert(vPSI,'set');
  sol:=pdsolve(syseqd, conjGerad, HINT=`+`);
  sol:=subs(map(i->i=op(0,i),conjGerad),sol); sol:=subs(PSI='Psi',sol);
  if nargs<6 or args[6]<>`all` then sol:=reduzConst(sol); fi;
  return sol;
end proc:
\end{verbatim}
\normalsize

\begin{description}
\item[Noether] given the infinitesimal generators which define a symmetry,
computes the conservation law for the optimal control problem, according
with Theorem~\ref{thm:TNoether} (Noether's theorem).
\item Output:
\begin{itemize}
\item[-] conservation law.
\end{itemize}
\item Syntax:
\begin{itemize}
\item[-] Noether(L, $\varphi$, t, x, u, S)
\end{itemize}
\item Input:
\begin{itemize}
\item[L -] expression of the Lagrangian;
\item[$\varphi$ -] expression or list of expressions of the velocity vector
$\boldsymbol{\varphi}$ which defines the control system;
\item[t -] name of the independent variable;
\item[x -] name or list of names of the state variables;
\item[u -] name or list of names of the control variables;
\item[S -] set of infinitesimal generators (output of procedure
\emph{Symmetry}).
\end{itemize}
\end{description}
\small
\begin{verbatim}
Noether:=proc(L::algebraic, phi::{algebraic, list(algebraic)}, t::name,
                            x0::{name,list(name)}, u0::{name,list(name)}, S::set)
  local n, xx, i, vX, vpsi, lpsi,  H, CL, lphi;
  unassign('T'); unassign('X'); unassign('psi');
  n:=nops(x0);
  if n>1 then lphi:=phi; lpsi:=[seq(psi[i],i=1..n)];
  else lpsi:=[psi]; lphi:=[phi]; fi;
  xx:=op(x0),op(u0),op(lpsi);
  vpsi:=Vector[row](lpsi);
  if n>1 then vX:=Vector([seq(X[i], i=1..n)]); else vX:=Vector([X]); fi;
  H:=psi_0*L+vpsi.Vector(lphi);
  CL:=vpsi.vX-H*T=const;
  CL:=eval(CL, S);
  CL:=subs({map(i->i=i(t),[xx])[]},CL); CL:=subs(psi_0=psi[0],CL);
  return CL;
end proc:
\end{verbatim}
\normalsize


\section{Concluding Remarks}
\label{concl}

Computer Algebra Systems are particularly suitable to handle the
problem of determining symmetries and conservation laws in optimal
control, because theory requires calculations that tend to be
tedious even for very simple problems with a linear control
system. \maple can perform these computations in a reliable way.
We illustrate our package by a large number of examples, that
range from simple problems in the calculus of variations to
optimal control problems of which the integrability was only very
recently shown by means of the conservation laws. Now these
computations can be made in a completely automatic way and without
any physical insights.

In mechanics, and calculus of variations, it is well known how to
use the conservation laws to reduce the order of the problems and,
with a sufficiently large number of independent conservation laws,
one can even integrate and solve the problems completely. Like
Noether's theorem, also the classical reduction theory can be
extended to the more general setting of optimal control
\cite{Echeverria,MR2005c:58032}. However, the reduction theory in
optimal control is an area not yet completed. More theoretical
results are needed in order to be possible to automatize the whole
process, from the computation of symmetries to the maximum
reduction of the problems.

As the examples show, the computing times
increase exponentially with the dimension
of the control system (with the number of state variables).
This is well illustrated with the problems
of sub-Riemannian geometry: problem $(2,3)$
(Example~\ref{ex:Ex10}), with three state variables,
requires a total computing time of one minute;
problem $(2,3,5)$ (Example~\ref{ex:Pr:235}),
with five state variables, has required
us thirty minutes. To tackle more complex problems a new approach is needed.
This is also under study and will be addressed elsewhere.
For example, integrability of the $(2,3,5,8)$ problem
of sub-Riemannian geometry, a problem with eight state variables,
is currently an open question. This problem
is out of the scope of the present
\maple package: eleven hours of computing time in our
1.4 GHz Pentium Centrino with 512MB of RAM
were not enough to determine the symmetries of the problem.


\section*{Acknowledgements}

PG was supported by the program PRODEP III/5.3/2003.
DT acknowledges the support from the control theory group
(\textsf{cotg}) of the R\&D unit \textsf{CEOC}, and the project
``Advances in Nonlinear Control and Calculus of Variations''
POCTI/MAT/41683/2001. The authors are grateful to E.~M.~Rocha
for stimulating discussions.



\begin{thebibliography}{99}

\bibitem{BaileyBorwein} D. H. Bailey and J. M. Borwein,
\textit{Experimental Mathematics: Examples, Methods and Implications},
Notices of the American Math. Society, {\bf 52} (2005), No.~5,
pp.~502--514.

\bibitem{Bonnard98} B. Bonnard, M. Chyba, and E. Tr\'{e}lat,
\textit{Sub-Riemannian Geometry: One-Parameter Deformation of the Martinet Flat Case},
Journal of Dynamical and Control Systems,
{\bf 4} (1998), No.~1, pp.~59--76.

\bibitem{BRU} B. van Brunt,
\textit{The Calculus of Variations},
Springer-Verlag New York, 2004.

\bibitem{Djukic73} D. S. Djukic,
\textit{Noether's theorem for optimum control systems},
Internat. J. Control, {\bf 1} (1973), No.~18, pp.~667--672.

\bibitem{Echeverria} A. Echeverr{\'{\i}}a-Enr{\'{\i}}quez,
J. Mar{\'{\i}}n-Solano, M. C. Mu\~{n}oz-Lecanda, and N. Rom\'{a}n-Roy,
\textit{Symmetries and reduction in optimal control theory},
Proceedings of the XI Fall Workshop on Geometry and Physics,
Publ. R. Soc. Mat. Esp., {\bf 6} (2004), pp.~203--208.

\bibitem{Gogodze88} I. K. Gogodze,
\textit{Symmetry in Problems of Optimal Control} (in Russian),
Proc. of extended sessions of seminar of the
Vekua Institute of Applied Mathematics,
Tbilisi University, Tbilisi,
{\bf 3} (1998), No.~3, pp.~39 -- 42.

\bibitem{gouv04}  P. D. F. Gouveia and D. F. M. Torres,
\textit{Computa\c{c}\~{a}o Alg\'{e}brica no C\'{a}lculo das Varia\c{c}\~{o}es:
Determina\c{c}\~{a}o de Simetrias e Leis de Conserva\c{c}\~{a}o} (in Portuguese),
Research report CM04/I-23, Dep. Mathematics, Univ. of Aveiro, September 2004.
Presented at XXVII CNMAC (Brazilian Congress of Applied Mathematics
and Computation), FAMAT/PUCRS, Porto Alegre, RS, Brasil, September 13-16, 2004.
\textsf{E-Print:} {\tt arXiv:math.OC/0411211}

\bibitem{GrizzleMarcus} J. W. Grizzle and S. I. Marcus,
\textit{The structure of nonlinear control systems possessing symmetries},
IEEE Trans. Automat. Control, {\bf 30} (1985), No.~3, pp.~248--258.

\bibitem{alik} A.~Gugushvili, O.~Khutsishvili, V.~Sesadze, G.~Dalakishvili,
N.~Mchedlishvili, T.~Khutsishvili, V.~Kekenadze, and D.~F.~M.~Torres,
\textit{Symmetries and Conservation Laws in Optimal Control Systems},
Georgian Technical University, Tbilisi, 2003.

\bibitem{LOG77} J. D. Logan,
\textit{Invariant Variational Principles},
Academic Press [Harcourt Brace Jovanovich Publishers], 1977.

\bibitem{LOG} J. D. Logan,
\textit{Applied Mathematics -- A Contemporary Approach},
John Wiley \& Sons, New York, 1987.

\bibitem{Rouchon01} Ph. Martin, R. M. Murray, and P. Rouchon,
\textit{Flat systems}, in {\it Mathematical control theory, Part 1, 2 (Trieste, 2001)},
ICTP Lect. Notes, VIII, Abdus Salam Int. Cent. Theoret. Phys., Trieste (2002),
pp.~705--768,

\bibitem{MR2005c:58032} E. Mart\'{\i}nez,
\textit{Reduction in optimal control theory},
Rep. Math. Phys. {\bf 53} (2004), No.~1, pp.~79--90.

\bibitem{Noether} E. Noether,
\textit{Invariante Variationsprobleme},
G\"{o}tt. Nachr. (1918), pp.~235--257.

\bibitem{OLV} P. J. Olver,
\textit{Applications of Lie Groups to Differential Equations},
Springer-Verlag, 1986.

\bibitem{IEEE05} B. Pal\'{a}ncz, Z. Beny\'{o}, and L. Kov\'{a}cs,
\textit{Control System Professional Suite, Product Review},
IEEE Control Systems Magazine {\bf 25} (2005), No.~4, pp.~67--75.

\bibitem{Pontryagin62} L. S. Pontryagin, V. G. Boltyanskii,
R. V. Gamkrelidze, and E. F. Mishchenko,
\textit{The mathematical theory of optimal processes},
Interscience Publishers John Wiley \& Sons, Inc. New York-London, 1962.

\bibitem{Richards} D. Richards,
\textit{Advanced Mathematical Methods with Maple},
Cambridge University Press, 2002.

\bibitem{Sachkov04} Yu. L. Sachkov,
\textit{Symmetries of flat rank two distributions and sub-Riemannian structures},
Trans. Amer. Math. Soc. {\bf 356} (2004), No.~2, pp.~457--494.

\bibitem{MR83c:70020} W. Sarlet and F. Cantrijn,
\textit{Generalizations of Noether's theorem in classical mechanics},
SIAM Rev., {\bf 23} (1981), No.~4, pp.~467--494.

\bibitem{Serovaiskii} S. Ya. Serovaiskii,
\textit{Counterexamples in optimal control theory},
Inverse and Ill-posed Problems Series, VSP, Utrecht, 2004.

\bibitem{Taimanov97} I. A. Taimanov,
\textit{Integrable geodesic flows of nonholonomic metrics},
J. Dynam. Control Systems, {\bf 3} (1997), No.~1, pp.~129--147.

\bibitem{delfimEJC} D. F. M. Torres,
\textit{On the Noether Theorem for Optimal Control},
European Journal of Control, {\bf 8} (2002), No.~1, pp.~56--63.

\bibitem{delfimIO} D. F. M. Torres,
\textit{A Remarkable Property of the Dynamic Optimization Extremals},
Investiga\c{c}\~{a}o Operacional, {\bf 22} (2002), No.~2, pp.~253--263.
\textsf{E-Print:} {\tt arXiv:math.OC/0212102}

\bibitem{3ncnw} D. F. M. Torres,
\textit{Conservation Laws in Optimal Control},
Dynamics, Bifurcations and Control,
F.~Colonius, L.~Gr\"{u}ne, eds.,
Lecture Notes in Control and Information Sciences {\bf 273} (2002),
Springer-Verlag, Berlin, Heidelberg, pp.~287--296.

\bibitem{TOR} D. F. M. Torres,
\textit{Proper Extensions of Noether's Symmetry Theorem
for Nonsmooth Extremals of the Calculus of Variations},
Communications on Pure and Applied Analysis,
{\bf 3} (2004), No.~3, pp.~491--500.
\textsf{E-Print:} {\tt arXiv:math.OC/0302127}

\bibitem{Torres04} D. F. M. Torres,
\textit{Quasi-Invariant Optimal Control Problems},
Portugali\ae\ Mathematica (N.S.), {\bf 61} (2004),
No.~1, pp.~97--114.
\textsf{E-Print:} {\tt arXiv:math.OC/0302264}

\bibitem{Torres05} D. F. M. Torres,
\textit{Weak Conservation Laws for Minimizers which are not
Pontryagin Extremals}, Proc. of the 2005 International Conference
``Physics and Control'' (PhysCon 2005), August 24-26, 2005, Saint
Petersburg, Russia.  Edited by
A.L.~Fradkov and A.N.~Churilov, 2005 IEEE, pp.~134--138.
\textsf{E-Print:} {\tt arXiv:math.OC/0503415}

\end{thebibliography}
\end{document}